\documentclass[11pt]{article}
\usepackage{amsmath,amsthm,amscd,amsfonts,amssymb,graphicx, color,mathrsfs}
\usepackage{mathtools}
\usepackage[colorlinks=true,linkcolor=blue,citecolor=blue,urlcolor=blue]{hyperref}
\usepackage{booktabs,tabularx,makecell, array}
\usepackage[utf8]{inputenc}
\usepackage{authblk} 
\usepackage{placeins}
\usepackage{float}
\usepackage[paper=a4paper, top=1in, bottom=1in, left=0.75in, right=0.75in]{geometry}
\newtheorem{theorem}{Theorem}[section]
\newtheorem{lemma}[theorem]{Lemma}
\theoremstyle{definition}
\newtheorem{definition}[theorem]{Definition}

\newtheorem{remark}[theorem]{Remark}

\numberwithin{equation}{section}
\numberwithin{table}{section}
\numberwithin{figure}{section}

\allowdisplaybreaks
\newcolumntype{L}[1]{>{\raggedright\arraybackslash}p{#1}} 
\newcolumntype{Y}{>{\raggedright\arraybackslash}X}        
\newcolumntype{M}{>{$\displaystyle}Y<{$}}

\begin{document}

\title{Normalized Solutions for Logarithmic Schr\"odinger Equations with Mass-Subcritical and HLS-Critical Choquard Nonlinearities}

\author[1]{Mohammad Saeid Abolhassanifar}
\author[1]{Reza Saadati}
\author[2]{Calogero Vetro\thanks{Corresponding Author}}
\author[1]{Mohammad Bagher Ghaemi}
\author[3]{Donal O'Regan}

\affil[1]{School of Mathematics and Computer Science, Iran University of Science and Technology, Narmak, Tehran 16846-13114, Iran}
\affil[2]{Department of Mathematics and Computer Science, University of Palermo, Via Archirafi 34, 90123 Palermo, Italy}
\affil[3]{School of Mathematical and Statistical Sciences, University of Galway, Galway, Ireland}

\date{}
\maketitle

\begin{center}
Emails: saied.aboul@gmail.com (M. S. Abolhassanifar); rsaadati@eml.cc (R. Saadati); calogero.vetro@unipa.it (C. Vetro); mghaemi@iust.ac.ir (M. B. Ghaemi); donal.oregan@nuigalway.ie (D.O'Regan).
\end{center}

\begin{abstract}
In this paper, we study the existence and multiplicity of normalized solutions to a Schrödinger equation involving a logarithmic nonlinearity, a Hardy–Littlewood–Sobolev critical Choquard nonlinearity, and a subcritical Choquard perturbation. To overcome the analytical difficulties caused by the interaction of these nonlinearities, we develop a new approach based on the geometry of the Pohozaev manifold. We obtain an explicit threshold \(c_0 > 0\) such that, for every \(0 < c < c_0\), {the degenerate part of the Pohozaev manifold is empty.} We first obtain a lower-energy solution by minimization on a gradient ball and then establish a second, higher-energy solution using an explicit energy gap. Both solutions are radial and strictly positive. A key feature of our energy gap analysis is that we do not need any restrictions on exponents, especially when $N \ge 5$. Our results correct and complement the non-existence conclusion in the mass-subcritical case considered in the previous paper [J. Math. Anal. Appl. 562 (2026), Paper No. 130707] and contribute to the study of logarithmic Schrödinger equations with multiple nonlinearities.
\end{abstract}



\section{Introduction}

\begin{sloppypar}
We consider normalized solutions of
\begin{align}\label{eq main}
\begin{cases}
-\Delta u + \lambda u = u\log|u|^2 + \mu(I_\alpha \ast |u|^p)|u|^{p - 2}u+(I_\beta \ast |u|^q)|u|^{q - 2}u & \text{in }\mathbb{R}^N,\\
\displaystyle\int_{\mathbb{R}^N}|u|^2\,dx = c^2,
\end{cases}
\end{align}
where $N\ge3$, $c>0$, $0<\alpha,\beta<N$, $\mu > 0$,
\[
\frac{N+\alpha}{N}<p<\frac{N+\alpha + 2}{N},
\qquad q=\frac{N+\beta}{N-2},
\]
and $\lambda\in\mathbb R$ is unknown. For $0<\sigma<N$, we define
\[
I_\sigma(x)=\frac{\Gamma(\frac{N-\sigma}{2})}{\Gamma(\frac{\sigma}{2})\pi^{N/2}2^\sigma}|x|^{\sigma-N}.
\]
Logarithmic Schr\"odinger equations originate in~\cite{Bialynicki1979}, while Choquard interactions arise in nonlocal quantum models; see~\cite{Moroz2017, Pekar1946}. Related normalized logarithmic and Choquard problems are treated in~\cite{Alves2020, BartschLiuLiu2020,   Cao2025, Deng2021, ShangMa2023, Shen2025, Shuai2025, XiaZhang2023, YaoChenRadulescuSun2022, Ye2022}.\\

\noindent
Under $t\star u(x)=e^{Nt/2}u(e^tx)$, the Dirichlet, lower Choquard, and critical Choquard terms scale as $e^{2t}$, $e^{(Np-N-\alpha)t}$, and $e^{2qt}$. Hence, the lower interaction is mass-subcritical whenever $Np-N-\alpha<2$. The $q$-term remains Hardy–Littlewood–Sobolev  (HLS, for short) critical. The logarithmic term is nonhomogeneous and is handled in the Orlicz setting of~\cite{Cazenave1983}.\\

\noindent
The closest logarithmic Choquard models are~\cite{MSA2026, Deng2025}. In~\cite{MSA2026}, both nonlocal terms use the same Riesz order. Here $I_\alpha$ and $I_\beta$ are independent: $Np-N-\alpha$ determines the mass geometry, whereas $I_\beta$ determines the critical compactness threshold.\\

\noindent
A step in the mass-subcritical argument of~\cite{MSA2026} requires correction. If $m:=Np-N-\alpha\in(0,2)$, substitution of $H'(t_0)=0$ into $H(t_0)$ gives the coefficient $2(q-1)/(2-m)>0$ in front of the critical Choquard term. Consequently, that computation alone does not imply $H(t_0)<0$; here the fiber signs are derived directly. By Lemma~\ref{lem 4.5}, we will see that every fiber has exactly two critical points. Compactness of the saddle branch follows from
\[
m^-(c)<m^+(c)+\frac{2+\beta}{2(N+\beta)}
S_\beta^{\frac{N+\beta}{2+\beta}},
\]
where $m^\pm(c)$ are admissible minimizer levels for the stable/unstable solutions (on the Pohozaev manifold). For $N=3,4$, the HLS optimizer is not in $L^2$, and we use a centered cut-off bubble; the decisive negative term is a mixed critical interaction. For $N\ge5$, the full HLS bubble lies in $L^2$ and is used without truncation. The logarithmic contribution $-C\varepsilon^2|\log\varepsilon|$ dominates the $O(\varepsilon^2)$ remainder, yielding a gap between $m^-(c)$ and $m^+(c)$ without additional exponent restrictions. The details are in the proof of Lemma~\ref{lem 4.10}. Related upper-critical normalized Choquard arguments appear in~\cite{ChenJinRadulescuZhou2026,ChenZhang2026,DingWang2024,GaoHe2024,LiuZhangSun2025,YangYang2025}.\\

\noindent
Our proof uses sharp nonlocal interpolation~\cite{Hongyu2016, Weinstein1982},
the HLS inequality~\cite{Lieb2001, Lieb1983}, and radial compactness in
the logarithmic--Orlicz space. Related variational and qualitative
Choquard results include~\cite{MaMoroz2024, MaMoroz2025, Moroz2013, WangMaLiu2023};
the positivity argument uses the strong maximum principle~\cite{Vazquez1984}.
Three features are decisive. First, the small-mass conditions give two
nondegenerate fibering critical points and a uniform bound $|Q|\ge\delta$
on the Pohozaev manifold. Second, the higher-energy solution is obtained
by minimizing directly on $\mathcal P_c^-$, using an Ekeland sequence
rather than a mountain-pass construction. Third, the strict gap
$m^-(c)<m^+(c)+\frac{q-1}{2q}S_\beta^{\frac q{q-1}}$ excludes an HLS bubble.
The mixed nonlocal interaction is decisive for $N=3,4$, while the
logarithmic term supplies the leading negative correction for $N\ge5$.

\end{sloppypar}

\subsection{Orlicz framework and functional setting}

We write $\|u\|_i$ for the $L^i(\mathbb{R}^N)$ norm of $u$ and $\|\nabla u\|_2$ for the $L^2$ norm of its gradient. Following \cite{Cazenave1983}, we define the Orlicz space
\[
G:= \left\{u \in L^1_{\mathrm{loc}}(\mathbb{R}^N) \mid \mathscr{A}(|u|)\in L^1(\mathbb{R}^N)\right\},
\] 
where the Young function $\mathscr{A}: [0, +\infty) \to [0, +\infty)$ is defined by
\[
\mathscr{A}(t):=
\begin{cases}
-t^2\log t^2, & 0\leq t<e^{-3},\\
3t^2 + 4e^{-3}t - e^{-6}, & t\geq e^{-3}.
\end{cases}
\]
We equip $G$ with the Luxemburg norm
\[
\|u\|_G:= \inf\left\{m > 0 \mid \int_{\mathbb{R}^N} \mathscr{A}\!\left(\frac{|u|}{m}\right)\,dx \leq 1\right\}.
\]
Define $H:= H^1(\mathbb{R}^N) \cap G$ with the norm $\|u\|_H:= \|u\|_{H^1(\mathbb{R}^N)} + \|u\|_G$, and
\begin{align}\label{de: log}
\mathscr{B}(t):= t^2\log t^2 + \mathscr{A}(t).
\end{align}
The functions $\mathscr{A}$ and $\mathscr{B}$ are nonnegative, convex, and increasing. For every $r>2$, there is $K_r>0$ such that
\begin{align}\label{ineq}
\mathscr{B}(t) \leq K_r t^r, \qquad (t \ge 0).
\end{align}
The space $H$ is a reflexive Banach space. We also recall the following standard properties of $G$:

\begin{itemize}
\item[$(i)$] for every $u\in G$,
\[
\inf\{\|u\|_G, \|u\|_G^2\} \leq \int_{\mathbb{R}^N} \mathscr{A}(|u|)\,dx \leq\sup\{\|u\|_G, \|u\|_G^2\};
\]
\item[$(ii)$] if $u_n \to u$ almost everywhere in $\mathbb{R}^N$ and
\[
\int_{\mathbb{R}^N} \mathscr{A}(|u_n|)\,dx \to \int_{\mathbb{R}^N} \mathscr{A}(|u|)\,dx<+\infty,
\]
then $\|u_n - u\|_G \to 0$.
\end{itemize}

\noindent
The radial subspace of $H$ is considered as $H_r:= H \cap H_r^1(\mathbb{R}^N)$ where
\[
H_r^1(\mathbb{R}^N):= \left\{u \in H^1(\mathbb{R}^N) \mid u(x)=u(|x|)\right\}.
\]
For $c > 0$, define
\[
S_c := \left\{u \in H_r \mid \int_{\mathbb{R}^N}|u|^2\,dx = c^2\right\}.
\]
Let $B_{\sigma, r}[u]:= \displaystyle\int_{\mathbb{R}^N} (I_\sigma \ast |u|^r)|u|^r\,dx$. For $c > 0$, we consider the radial mass sphere $S_c$. The energy functional for problem \eqref{eq main} is defined as:
\begin{align}\label{energy}
E(u) &:= \frac{1}{2}\|\nabla u\|_2^2 + \frac{1}{2}\|u\|_2^2 - \frac{1}{2}\int_{\mathbb{R}^N} |u|^2\log{|u|^2}\,dx \notag \\
& \quad - \frac{\mu}{2p}B_{\alpha, p}[u] - \frac{1}{2q} B_{\beta, q}[u].
\end{align}

\begin{definition}\label{def 1}
A function $u \in S_c$ is a normalized critical point of $E$ in \eqref{energy} if $d(E|_{S_c})(u) = 0$; equivalently, there exists $\lambda\in\mathbb R$ such that $E'(u)+\lambda\Psi'(u)=0$, where $\Psi(u):=\frac12\|u\|_2^2$.\\

\noindent
A normalized critical point is called a radial normalized ground state if
\[
E(u^\ast) = \inf\left\{E(u) \mid u \in S_c,\, d(E|_{S_c})(u) = 0\right\},
\]
and a radially excited state if
\[
E(u^\ast) > \inf\left\{E(u) \mid u \in S_c,\, d(E|_{S_c})(u) = 0\right\},
\]
where the first variation of $E$ is
\begin{align}\label{energy: variation}
E'(u)[v] &=\int_{\mathbb{R}^N} \nabla u\cdot \nabla v\,dx - \int_{\mathbb{R}^N} u\log|u|^2\,v\,dx \notag\\
&\quad-\mu\int_{\mathbb{R}^N} (I_\alpha \ast |u|^p)|u|^{p - 2}u\,v\,dx - \int_{\mathbb{R}^N} (I_\beta \ast |u|^q)|u|^{q - 2}u\,v\,dx,
\end{align}
for every $v\in H$. Since the functional and the constraint are invariant under the orthogonal group, the principle of symmetric criticality implies that every radial normalized critical point is a weak solution of \eqref{eq main} in the full test space.
\end{definition}

\noindent
To analyze existence and multiplicity, we introduce the Pohozaev manifold
\[
\mathcal{P}_c := \left\{u \in S_c \mid P(u) = 0 \right\},
\]
with
\begin{align}\label{P}
P(u) &:= \|\nabla u\|_2^2 - \frac{N c^2}{2} - \mu\eta_pB_{\alpha, p}[u] - B_{\beta, q}[u],
\end{align}
where $\eta_p:= \frac{Np - N - \alpha}{2p}$. For $x \in \mathbb{R}^N$, we define the mass preserving scale by $t \star u(x):= e^{\frac{N}{2}t} u(e^t x)$. According to this scaling, we consider the fibering map associated to $E$ in \eqref{energy} by
\begin{align}\label{fib}
\Phi_u(t) &:= E(t \star u) \notag \\
& = \frac{e^{2t}}{2}\|\nabla u\|_2^2 + \frac{1 - Nt}{2}c^2  - \frac{1}{2}\int_{\mathbb{R}^N} |u|^2 \log{|u|^2}\,dx\notag \\
& \quad - \frac{\mu e^{2p\eta_p t}}{2p}B_{\alpha, p}[u] - \frac{e^{2q t}}{2q}B_{\beta, q}[u],
\end{align}
where $\Phi_u'(t) = P(t\star u)$. To analyze the multiplicity and stability of solutions, and to determine which is the ground state or excited state, it is useful to split $\mathcal{P}_c$ into its three components as $\mathcal{P}_c = \mathcal{P}_c^+ \cup \mathcal{P}_c^0 \cup \mathcal{P}_c^-$, with
\[
\begin{aligned}
\mathcal{P}_c^+ &:= \left\{u \in S_c \mid \Phi_u'(0)=0,\, \Phi_u''(0) > 0 \right\}, \\
\mathcal{P}_c^{0} &:= \left\{ u \in S_c \mid \Phi_u'(0)=0,\, \Phi_u''(0) = 0 \right\}, \\
\mathcal{P}_c^{-} &:= \left\{ u \in S_c \mid \Phi_u'(0)=0,\, \Phi_u''(0) < 0 \right\},
\end{aligned}
\]
where
\begin{align}\label{Q}
Q(u):= \Phi_u''(0) = 2\|\nabla u\|_2^2 - 2\mu p\eta_p^2B_{\alpha, p}[u] - 2qB_{\beta, q}[u].
\end{align}

\subsection{Main results}

We now state the main results. We define the ground variational level and the levels associated with the two Poho\v{z}aev branches by
\[
m(c):=\inf_{u\in\mathcal{P}_c}E(u),\qquad
m^+(c):= \inf_{u \in \mathcal{P}_c^+ \cap H_r} E(u), \qquad
m^-(c):= \inf_{u \in \mathcal{P}_c^- \cap H_r} E(u).
\]
Also, due to the existence of the HLS critical term, we need to consider the following constant
\begin{align}\label{S}
S_\beta:= \inf_{u \in D^{1, 2}(\mathbb{R}^N) \setminus \{0\}} \frac{\|\nabla u\|_2^2}{\left(\displaystyle\int_{\mathbb{R}^N}(I_\beta \ast|u|^q)|u|^q\,dx\right)^{\frac{1}{q}}}.
\end{align}
This sharp constant determines the one-bubble threshold. Set $c_0:=\min\{c_F,c_G,c_K\}>0$, with the thresholds specified in Section~\ref{sec:existence}.

\begin{theorem}\label{thm 1}
Let $N\ge3$, $0<\alpha,\beta<N$, and
\[
\frac{N+\alpha}{N}<p<\frac{N+\alpha+2}{N},
\qquad q=\frac{N+\beta}{N-2}.
\]
Assume that $\mu>0$, $K>\frac{q}{q-1}$, and
$0<c<c_0$. Then problem~\eqref{eq main}
admits at least two distinct positive radial normalized solutions
$u_c^+,u_c^-\in S_c$, with Lagrange multipliers
$\lambda_c^+,\lambda_c^-\in\mathbb R$, satisfying
\[
\begin{gathered}
 u_c^+\in\mathcal P_c^+,
 \qquad E(u_c^+)=m^+(c),\qquad Q(u_c^+)>0,\\
 u_c^-\in\mathcal P_c^-,
 \qquad E(u_c^-)=m^-(c),\qquad Q(u_c^-)<0.
\end{gathered}
\]
Moreover,
\[
 m^+(c)<m^-(c)
 <m^+(c)+\frac{\beta+2}{2(N+\beta)}
 S_\beta^{\frac{N+\beta}{\beta+2}}.
\]
\end{theorem}

\section{Preliminaries}\label{sec 2}

In this section, we recall some tools that are needed for our analysis. We begin with the $L^2$ compactness supplied by the radial logarithmic--Orlicz setting.
\begin{lemma}\label{lem 2.1}
The embedding $H_r\hookrightarrow L^2(\mathbb R^N)$ is compact.
\end{lemma}

\begin{proof}
Let $(u_n)$ be a bounded sequence in $H_r$ and, after passing to a subsequence, assume that $u_n\rightharpoonup u$ in $H_r$. By the radial lemma, we have
\[
|u_n(x)|^2 \leq C|x|^{1 - N}\|u_n\|_2\|\nabla u_n\|_2 \leq C_0|x|^{1 - N},
\]
uniformly in $n$ for $x\neq 0$. Hence, for every sufficiently large $R$, we have $|u_n(x)| < e^{-3}$ on $|x| > R$. Hence,
\[
\mathscr{A}(|u_n(x)|) = -|u_n(x)|^2\log |u_n(x)|^2 \geq c \log R\,|u_n(x)|^2,
\]
where $c > 0$ is independent of $n$ and $R$ for $|x|>R$. By boundedness in $G$ and property~$(i)$, we obtain $\sup_{n \ge 0} \int \mathscr{A}(|u_n|) < +\infty$ and
\[
\sup_{n \ge 0} \int_{|x| > R}|u_n|^2\,dx \leq \frac{C}{\log R} \to 0,
\]
as $R \to +\infty$. The same tail estimate holds for $u$. On each fixed ball $B_R$, Rellich's compact embedding theorem gives $u_n\to u$ in $L^2(B_R)$ after extraction. Combining local convergence with the uniform tail estimate yields $u_n \to u$ in $L^2(\mathbb{R}^N)$.
\end{proof}

\noindent
We use the Gagliardo--Nirenberg inequality~\cite{Weinstein1982}: for $u\in H^1(\mathbb{R}^N)$ and $2 < r < 2^*:=\frac{2N}{N - 2}$,
\begin{align}\label{GN}
\|u\|_r^r \leq C_{N, r} \|\nabla u\|_2^{r\gamma_r}\,\|u\|_2^{r(1 - \gamma_r)}.
\end{align}
Here $\gamma_r := N\left(\frac{1}{2} - \frac{1}{r}\right)$ and $C_{N, r} > 0$ is the sharp constant. We also use the HLS inequality~\cite{Lieb2001}. If $N \ge 1$, $s, t > 1$, $0 < \sigma < N$, and $\frac{1}{s} + \frac{N - \sigma}{N} + \frac{1}{t} = 2$, then
\[
\left|\int_{\mathbb{R}^N}\int_{\mathbb{R}^N}\frac{u(x)v(y)}{|x - y|^{N - \sigma}}\,dx\,dy\right| \leq C(N, \sigma, s)\|u\|_s \|v\|_t,
\]
for $u \in L^s(\mathbb{R}^N)$ and $v \in L^t(\mathbb{R}^N)$.\\

\noindent
Let $N \ge 3$, $0 < \sigma < N$, and $\frac{N + \sigma}{N} < r < \frac{N + \sigma}{N - 2}$. We consider
\begin{align}\label{Critical: Form}
-\frac{Nr - (N + \sigma)}{2}\Delta Q_{\sigma, r} +\frac{N + \sigma - (N - 2)r}{2}Q_{\sigma, r} = \left(I_\sigma \ast|Q_{\sigma, r}|^r\right)|Q_{\sigma, r}|^{r - 2}Q_{\sigma, r}.
\end{align}
For $N \ge 3$, $0 < \sigma < N$, and $\frac{N + \sigma}{N} < r < \frac{N + \sigma}{N - 2}$, the sharp nonlocal version of the Gagliardo--Nirenberg inequality (see \cite{Hongyu2016}) is
\begin{align}\label{GS}
B_{\sigma, r}[u] \leq C_{N, \sigma, r} \|\nabla u\|_2^{2r\eta_r}\|u\|_2^{2r - 2r\eta_r},
\end{align}
where $\eta_r=\frac{Nr-N-\sigma}{2r}$ and $C_{N, \sigma, r} = \frac{r}{\|Q_{\sigma, r}\|_2^{2(r - 1)}}$ and $Q_{\sigma, r}$ is the positive optimizer and a weak solution of \eqref{Critical: Form}.

\begin{lemma}\label{lem 2.2}\cite{Moroz2013}
Let $N \ge 3$, $0 < \sigma < N$, and let $r$ belong to the admissible HLS range. Let $(v_n) \subset H^1(\mathbb{R}^N)$ be a bounded sequence satisfying
$v_n \rightharpoonup v$ weakly in $H^1(\mathbb{R}^N)$, and $v_n \to v$ almost everywhere in $\mathbb{R}^N$. Set $w_n = v_n - v$, then
\begin{align}\label{Choquard: convergence}
B_{\sigma, r}[v_n] - B_{\sigma, r}[w_n] \to B_{\sigma, r}[v].
\end{align}
\end{lemma}

\begin{lemma}\label{lem 2.3}\cite{Shen2025}
Let $p \ge 2$ be an exponent. Let $(v_n) \subset H$ be a bounded sequence satisfying $v_n \to v$ almost everywhere in $\mathbb{R}^N$, and the sequence $|v_n|^p \log|v_n|^p \in L^1(\mathbb{R}^N)$ be bounded. Set $w_n = v_n - v$, then
\begin{align}\label{Log: convergence}
 \int_{\mathbb{R}^N}|v_n|^p \log|v_n|^p\,dx - \int_{\mathbb{R}^N}|w_n|^p \log|w_n|^p\,dx \to \int_{\mathbb{R}^N}|v|^p \log|v|^p\,dx.
\end{align}
\end{lemma}

\noindent
Differentiating \eqref{P} for $u,v\in H_r$, we obtain
\begin{align*}
P'(u)[v] &= 2\int_{\mathbb{R}^N}\nabla u \cdot \nabla v\,dx - 2\mu p\eta_p\int_{\mathbb{R}^N}(I_\alpha \ast |u|^p)|u|^{p - 2}uv\,dx \notag \\
&\quad -2q\int_{\mathbb{R}^N}(I_\beta \ast |u|^q)|u|^{q - 2}uv\,dx.
\end{align*}
By the HLS and Sobolev inequalities, we obtain
\begin{align*}
\|P'(u)\|_{H_r^*} &\leq C\left(\|u\|_{H_r} + \|u\|_{H_r}^{2p - 1} + \|u\|_{H_r}^{2q - 1}\right).
\end{align*}
Thus $P'(u_n)$ is bounded in $H_r^*$ whenever $(u_n)$ is bounded in $H_r$.\\

\noindent
Regarding the convergence of the Young function $\mathscr{A}$, for every $u \in S_c$, we can write
\begin{align}\label{orlicz-recovery-identity}
\int_{\mathbb R^N}\mathscr{A}(|u|)\,dx &= 2E(u) - \|\nabla u\|_2^2 - c^2 + \int_{\mathbb R^N}\mathscr{B}(|u|)\,dx \notag \\
&\quad + \frac{\mu}{p}B_{\alpha, p}[u] +\frac{1}{q}B_{\beta, q}[u].
\end{align}
If $u_n \to u$ in $H^1(\mathbb{R}^N)$, $E(u_n) \to E(u)$, and the $B$-term and both Choquard terms converge, then
\[
\int_{\mathbb{R}^N}\mathscr{A}(|u_n|)\,dx \to \int_{\mathbb{R}^N}\mathscr{A}(|u|)\,dx,
\]
as $n \to +\infty$.

\begin{lemma}\label{lem 2.5}
Let $N \ge 3$, $0 < \alpha < N$, and $\frac{N + \alpha}{N} < p < \frac{N + \alpha}{N - 2}$. Let $(u_n) \subset H_r$ be uniformly bounded and satisfy $u_n \rightharpoonup u$ in $H_r$ and $u_n \to u$ in $L^2(\mathbb{R}^N)$. Then
\begin{align}\label{subcrit-compact-p}
B_{\alpha, p}[u_n] \to B_{\alpha, p}[u],
\end{align}
as $n \to +\infty$. Moreover, for every $v \in H_r$,
\begin{align}\label{subcrit-compact-p-variation}
\int_{\mathbb{R}^N}(I_\alpha \ast |u_n|^p)|u_n|^{p - 2}u_n v\,dx \to \int_{\mathbb{R}^N}(I_\alpha \ast |u|^p)|u|^{p - 2}u v\,dx,
\end{align}
as $n \to +\infty$.
\end{lemma}

\begin{proof}
By $2 < \frac{2Np}{N + \alpha} < 2^*$, the radial compact embedding and interpolation, we obtain $u_n \to u$ in $L^{\frac{2Np}{N + \alpha}}(\mathbb{R}^N)$. Hence $|u_n|^p \to |u|^p$ in $L^{\frac{2N}{N + \alpha}}(\mathbb{R}^N)$, and applying the HLS inequality, we obtain
\begin{align*}
\left|B_{\alpha, p}[u_n] - B_{\alpha, p}[u]\right| \leq C \big\||u_n|^p - |u|^p\big\|_{\frac{2N}{N + \alpha}}\left(\big\||u_n|^p\big\|_{\frac{2N}{N + \alpha}} + \big\||u|^p\big\|_{\frac{2N}{N + \alpha}}\right) \to 0,
\end{align*}
as $n \to +\infty$. The Nemytskii map $z \mapsto |z|^{p - 2}z$ is continuous for every $p > 1$. Hence $|u_n|^{p - 2}u_n \to |u|^{p - 2}u$ in $L^{\frac{2Np}{(N+\alpha)(p - 1)}}(\mathbb{R}^N)$, as $n \to +\infty$. Again, we apply the HLS inequality, and we obtain
\begin{align*}
&\left| \int_{\mathbb{R}^N}(I_\alpha \ast |u_n|^p)|u_n|^{p - 2}u_n v\,dx - \int_{\mathbb R^N}(I_\alpha \ast |u|^p)|u|^{p - 2}u v\,dx\right| \\
&\leq C \big\||u_n|^p - |u|^p\big\|_{\frac{2N}{N + \alpha}}\big\||u_n|^{p - 2}u_n\big\|_{\frac{2Np}{(N + \alpha)(p - 1)}}\|v\|_{\frac{2Np}{N + \alpha}}\\
&\quad + C\big\||u|^p\big\|_{\frac{2N}{N + \alpha}}\big\||u_n|^{p - 2}u_n - |u|^{p - 2}u\big\|_{\frac{2Np}{(N+\alpha)(p - 1)}}\|v\|_{\frac{2Np}{N + \alpha}} \to 0,
\end{align*}
as $n \to +\infty$. This completes the proof.
\end{proof}

\begin{lemma}\label{lem 2.6}
Let $u\in H$ be a weak solution of \eqref{eq main}.
Then $u\in L^\infty_{\mathrm{loc}}(\mathbb R^N)$.
\end{lemma}

\begin{proof}
For $(\gamma,r)=(\alpha,p),(\beta,q)$, the range $\frac{N+\gamma}{N}<r\le\frac{N+\gamma}{N-2}$ and $u\in H^1(\mathbb R^N)$ give, for $A\ge1$,
\[
\begin{gathered}
|u|^{r-1}\mathbf1_{\{|u|\le A\}}\in L^{\frac{2N}{\gamma}},
\qquad
|u|^{r-1}\mathbf1_{\{|u|>A\}}\in L^{\frac{2N}{\gamma+2}},\\
\bigl\||u|^{r-1}\mathbf1_{\{|u|>A\}}\bigr\|_{\frac{2N}{\gamma+2}}
\to0\quad(A\to\infty),
\qquad
(\log |u|^2)_+\in L^{N/2}.
\end{gathered}
\]
For $z\in H^1(\mathbb R^N)$ and $0<\theta\le1$ close to $1$, HLS, H\"older, Sobolev inequalities and interpolation between $L^2$ and $L^{2^*}$ give
\[
\begin{aligned}
&\int_{\mathbb R^N}
 \bigl(I_\gamma*(|u|^{r-1}|z|^\theta)\bigr)
 |u|^{r-1}|z|^{2-\theta}\,dx\\
&\quad\le
 C\bigl\||u|^{r-1}\mathbf1_{\{|u|>A\}}\bigr\|_{\frac{2N}{\gamma+2}}^2
 \|\nabla z\|_2^2\\
&\qquad+
 C\bigl\||u|^{r-1}\mathbf1_{\{|u|>A\}}\bigr\|_{\frac{2N}{\gamma+2}}
 \bigl\||u|^{r-1}\mathbf1_{\{|u|\le A\}}\bigr\|_{\frac{2N}{\gamma}}
 \|\nabla z\|_2\|z\|_2\\
&\qquad+
 C\bigl\||u|^{r-1}\mathbf1_{\{|u|\le A\}}\bigr\|_{\frac{2N}{\gamma}}^2
 \|z\|_2^2.
\end{aligned}
\]
Since
\[
\int_{\mathbb R^N}(\log |u|^2)_+z^2\,dx
\le\varepsilon\|\nabla z\|_2^2+C_\varepsilon\|z\|_2^2,
\]
choosing $A$ large and applying Young's inequality yields
\[
\begin{aligned}
&\int_{\mathbb R^N}(\log |u|^2)_+z^2\,dx\\
&\quad+\mu\int_{\mathbb R^N}
 \bigl(I_\alpha*(|u|^{p-1}|z|^\theta)\bigr)
 |u|^{p-1}|z|^{2-\theta}\,dx\\
&\quad+\int_{\mathbb R^N}
 \bigl(I_\beta*(|u|^{q-1}|z|^\theta)\bigr)
 |u|^{q-1}|z|^{2-\theta}\,dx\\
&\le\varepsilon\|\nabla z\|_2^2+C_{\varepsilon,\theta}\|z\|_2^2.
\end{aligned}
\]
The same bound holds for pointwise truncations of $(\log |u|^2)_+$, $\sqrt{\mu}|u|^{p-1}$ and $|u|^{q-1}$, with the same constants. Set $T_k(f)=\min\{f,k\}\mathbf1_{B_k}$ for $f\ge0$. For $\Lambda>\max\{1,\lambda\}$ sufficiently large, consider
\[
\begin{aligned}
-\Delta v_k+\Lambda v_k={}&T_k((\log |u|^2)_+)v_k\\
&+T_k(\sqrt{\mu}|u|^{p-1})
 \bigl(I_\alpha*(T_k(\sqrt{\mu}|u|^{p-1})v_k)\bigr)\\
&+T_k(|u|^{q-1})
 \bigl(I_\beta*(T_k(|u|^{q-1})v_k)\bigr)
 +(\Lambda-\lambda)|u|.
\end{aligned}
\]
The form bound with $\theta=1$ gives
\[
\begin{aligned}
&\|\nabla z\|_2^2+\Lambda\|z\|_2^2
 -\int_{\mathbb R^N}T_k((\log |u|^2)_+)z^2\,dx\\
&\quad-\int_{\mathbb R^N}
 \bigl(I_\alpha*(T_k(\sqrt{\mu}|u|^{p-1})z)\bigr)
 T_k(\sqrt{\mu}|u|^{p-1})z\,dx\\
&\quad-\int_{\mathbb R^N}
 \bigl(I_\beta*(T_k(|u|^{q-1})z)\bigr)
 T_k(|u|^{q-1})z\,dx\\
&\ge\tfrac12\|\nabla z\|_2^2+\|z\|_2^2.
\end{aligned}
\]
By the Lax--Milgram theorem, we have a unique $v_k\in H^1(\mathbb R^N)$. Testing with $v_k^-$ and then $v_k$ yields $v_k\ge0$ and $\sup_k\|v_k\|_{H^1}<\infty$. Choose $2<m<\min\left\{2^*,\frac{2N}{\alpha},\frac{2N}{\beta}\right\}$ close enough to $2$ so that $\theta=2/m$ is admissible. Then $\sup_k\|v_k\|_m<\infty$. Also, for $v_{k,L}=\min\{v_k,L\}$, testing with $v_{k,L}^{m-1}$ gives
\[
\begin{aligned}
&\frac{4(m-1)}{m^2}\|\nabla v_{k,L}^{m/2}\|_2^2\\
&\quad\le\int_{\mathbb R^N}T_k((\log |u|^2)_+)v_{k,L}^m\,dx\\
&\qquad+\int_{\mathbb R^N}
 \bigl(I_\alpha*(T_k(\sqrt{\mu}|u|^{p-1})v_{k,L})\bigr)
 T_k(\sqrt{\mu}|u|^{p-1})v_{k,L}^{m-1}\,dx\\
&\qquad+\int_{\mathbb R^N}
 \bigl(I_\beta*(T_k(|u|^{q-1})v_{k,L})\bigr)
 T_k(|u|^{q-1})v_{k,L}^{m-1}\,dx\\
&\qquad+(\Lambda-\lambda)
 \int_{\mathbb R^N}|u|v_{k,L}^{m-1}\,dx+R_{k,L}.
\end{aligned}
\]
By positivity and symmetry of the kernel, we get
\[
\begin{aligned}
0\le R_{k,L}\le{}&
 \int_{\{v_k>L\}}T_k((\log |u|^2)_+)v_k^m\,dx\\
&+\int_{\{v_k>L\}}
 \bigl(I_\alpha*(T_k(\sqrt{\mu}|u|^{p-1})v_k^{m-1})\bigr)
 T_k(\sqrt{\mu}|u|^{p-1})v_k\,dx\\
&+\int_{\{v_k>L\}}
 \bigl(I_\beta*(T_k(|u|^{q-1})v_k^{m-1})\bigr)
 T_k(|u|^{q-1})v_k\,dx.
\end{aligned}
\]
By HLS inequality with reciprocal exponents
$\frac{\gamma}{2N}+\frac1m$ and
$\frac{\gamma}{2N}+1-\frac1m$, and by H\"older's inequality, for fixed $k$, we obtain
\[
R_{k,L}\le
 k\|v_k\mathbf1_{\{v_k>L\}}\|_m^m
 +C_k\|v_k\|_m^{m-1}\|v_k\mathbf1_{\{v_k>L\}}\|_m
 \to0\qquad(L\to\infty).
\]
The form bound with $z=v_{k,L}^{m/2}$ and $\theta=2/m$ yields
\[
\frac{2(m-1)}{m^2}\|\nabla v_{k,L}^{m/2}\|_2^2
\le C_m\bigl(\|v_k\|_m^m+\|u\|_m^m\bigr)+R_{k,L}.
\]
By Sobolev inequality and $L\to\infty$ for each $k$, we have $\sup_k\|v_k\|_{\frac{mN}{N-2}}<\infty$. The coefficient convergences
\[
\begin{cases}
T_k((\log |u|^2)_+)\to(\log |u|^2)_+ & \text{in } L^{N/2},\\
T_k(\sqrt{\mu}|u|^{p-1})\to\sqrt{\mu}|u|^{p-1}&\text{in } L^{\frac{2N}{\alpha+2}}+L^{\frac{2N}{\alpha}},\\
T_k(|u|^{q-1})\to|u|^{q-1} &\text{in } L^{\frac{2N}{\beta+2}}+L^{\frac{2N}{\beta}},
\end{cases}
\]
allow passage to the limit along a subsequence $v_k\rightharpoonup v$ in $H^1$, with $v\in L^{\frac{mN}{N-2}}$. We then obtain
\[
\begin{aligned}
-\Delta v+\Lambda v={}&(\log |u|^2)_+v\\
&+\mu |u|^{p-1}\bigl(I_\alpha*(|u|^{p-1}v)\bigr)
 +|u|^{q-1}\bigl(I_\beta*(|u|^{q-1}v)\bigr)
 +(\Lambda-\lambda)|u|.
\end{aligned}
\]
By Kato's inequality, $|u|$ is a subsolution of this equation.
Testing the difference with $(|u|-v)_+$ and using kernel positivity
and coercivity gives
\[
\tfrac12\|\nabla(|u|-v)_+\|_2^2+\|(|u|-v)_+\|_2^2\le0.
\]
Thus $|u|\le v$ and $u\in L^{\frac{mN}{N-2}}(\mathbb R^N)$,
where $\frac{mN}{N-2}>2^*$. Fix $B_\rho(x_0)\Subset B_R(x_0)$. For each $(\gamma,r)$,
\[
\sup_{x\in B_\rho(x_0)}
 \int_{\mathbb R^N\setminus B_R(x_0)}I_\gamma(x-y)|u(y)|^r\,dy
 \le C_{R-\rho}\|u\|_{\frac{2Nr}{N+\gamma}}^r<\infty.
\]
Suppose $u\in L^s(B_R(x_0))$ with $s\ge\frac{mN}{N-2}$.
For $s<Nr/\gamma$, the HLS inequality leads to
\[
\left\|I_\gamma*(|u|^r\mathbf1_{B_R(x_0)})\right\|
 _{L^{\frac{Ns}{Nr-\gamma s}}(B_\rho(x_0))}
 \le C\|u\|_{L^s(B_R(x_0))}^r.
\]
At $s=Nr/\gamma$ the interior potential belongs to every finite
$L^a(B_\rho(x_0))$, and for $s>Nr/\gamma$ it is bounded there.
Hence $(I_\gamma*|u|^r)|u|^{r-2}u\in L^b(B_\rho(x_0))$ for every $b>1$ with
\[
\frac1b>
\left(\frac rs-\frac\gamma N\right)_++\frac{r-1}s.
\]
Also, $|t\log t^2|\le C_\delta(1+|t|^{1+\delta})$ gives
$u\log |u|^2\in L^{s/(1+\delta)}(B_\rho(x_0))$. Now
\[
\begin{aligned}
&\frac1s+\frac2N-
 \left[\left(\frac rs-\frac\gamma N\right)_++\frac{r-1}s\right]
=
\begin{cases}
\displaystyle\frac{\gamma+2}{N}-\frac{2(r-1)}s,
 &s\le Nr/\gamma,\\[2mm]
\displaystyle\frac2N-\frac{r-2}s,
 &s\ge Nr/\gamma.
\end{cases}
\end{aligned}
\]
Both expressions are uniformly positive for $s\ge\frac{mN}{N-2}$
because $r-1\le\frac{\gamma+2}{N-2}$, $\frac{mN}{N-2}>2^*$
and $r<2^*$. For small fixed $\delta>0$, the logarithmic term
gives the uniformly positive gap $\frac2N-\frac{\delta}{s}$.\\

\noindent
Local $W^{2,b}$ estimates and Sobolev embedding therefore decrease
$1/s$ by a fixed positive amount on nested balls. After finitely
many steps, the right-hand side of \eqref{eq main} belongs to
$L^b_{\mathrm{loc}}$ for some $b>N/2$. Hence
$u\in L^\infty_{\mathrm{loc}}(\mathbb R^N)$.
\end{proof}

\begin{remark}\label{rem 2.6}
By Lemma~\ref{lem 2.6}, we have $u\in L^\infty_{\mathrm{loc}}(\mathbb R^N)$. {Since the Riesz convolution terms considered in \eqref{eq main} are locally bounded, the right-hand side of \eqref{eq main} is bounded too}, and $-\Delta u\in L^\infty_{\mathrm{loc}}(\mathbb R^N)$. Hence, by local elliptic regularity and Sobolev embedding, we have  $u\in C^{1,\theta}_{\mathrm{loc}}(\mathbb R^N)$ for every $\theta \in (0, 1)$.
\end{remark}

{
\begin{lemma}\label{lem 2.4}
Let $u\in S_c$. Then the following assertions hold.
\begin{enumerate}
\item[$(i)$] If $u$ is a normalized critical point of $E$ on $S_c$, then $P(u)=0$.
\item[$(ii)$] If $u\in\mathcal P_c$ is a critical point of $E|_{\mathcal P_c}$
and $Q(u)=\Phi_u''(0)\ne0$, then $u$ is a normalized
critical point of $E$ on $S_c$.
\end{enumerate}
\end{lemma}

\begin{proof}
\noindent\textit{Part $(i)$.}
By symmetric criticality, $u$ solves \eqref{eq main} weakly for some
$\lambda\in\mathbb R$. Lemma~\ref{lem 2.6} and Remark~\ref{rem 2.6} yield
\[
 u\in W^{2,s}_{\mathrm{loc}}(\mathbb R^N)
       \cap C^{1,\theta}_{\mathrm{loc}}(\mathbb R^N)
 \qquad (1<s<\infty,\ 0<\theta<1).
\]
Since $u\in H_r$ need not imply $x\cdot\nabla u\in H^1$
(see \cite[Remark~2.5]{BartschSoave2017}), take radial
$\chi\in C_c^\infty$, $0\le\chi\le1$, equal to $1$ on $B_1$
and $0$ outside $B_2$. For $R>1$, set
\[
 \chi_R(x)=\chi(x/R),\qquad V_R(x)=x\chi_R(x),\qquad
 Z_R(x)=V_R(x)\cdot\nabla u(x).
\]
Local regularity and $\operatorname{supp}Z_R\subset B_{2R}$ give
\[
\begin{aligned}
 \|Z_R\|_2
 &\le 2R\|\nabla u\|_{L^2(B_{2R})}<\infty,\\
 \|\nabla Z_R\|_2
 &\le \|\nabla V_R\|_\infty\|\nabla u\|_{L^2(B_{2R})}
       +2R\|D^2u\|_{L^2(B_{2R})},\\
 &\le C\|\nabla u\|_{L^2(B_{2R})}
       +2R\|D^2u\|_{L^2(B_{2R})}<\infty.
\end{aligned}
\]
Since $Z_R$ is radial, bounded, and compactly supported,
$\int_{\mathbb R^N}\mathscr A(|Z_R|)\,dx<\infty$; thus $Z_R\in H_r$ for each $R$.
No uniform bound is needed. Testing with $u$ gives
\begin{equation}\label{eq:lemma24-test-u}
\begin{aligned}
 \|\nabla u\|_2^2+\lambda c^2
 &=\int_{\mathbb R^N}|u|^2\log|u|^2\,dx
       +\mu B_{\alpha,p}[u]+B_{\beta,q}[u],\\
 \int_{\mathbb R^N}|u|^2\log|u|^2\,dx
 &=\|\nabla u\|_2^2+\lambda c^2
       -\mu B_{\alpha,p}[u]-B_{\beta,q}[u].
\end{aligned}
\end{equation}
By direct calculation with $Z_R$, as $R\to\infty$, we obtain
\[
\begin{aligned}
 \int_{\mathbb R^N}(-\Delta u)Z_R\,dx
 &\longrightarrow \frac{2-N}{2}\|\nabla u\|_2^2,\\
 \int_{\mathbb R^N}(\lambda u-u\log|u|^2)Z_R\,dx
 &\longrightarrow -\frac{\lambda N}{2}c^2
  +\frac N2\int_{\mathbb R^N}|u|^2\log|u|^2\,dx
  -\frac N2c^2,\\
 \int_{\mathbb R^N}(I_\sigma*|u|^r)|u|^{r-2}u\,Z_R\,dx
 &\longrightarrow -\frac{N+\sigma}{2r}B_{\sigma,r}[u],
 \qquad (\sigma,r)\in\{(\alpha,p),(\beta,q)\}.
\end{aligned}
\]
The cutoff errors vanish by integrability; the nonlocal limits use
HLS, symmetrization, and dominated convergence as in
\cite[Proposition~3.1]{Moroz2013}.
Combining these limits with \eqref{eq:lemma24-test-u} gives
\[
\begin{aligned}
0={}&\frac{2-N}{2}\|\nabla u\|_2^2
      -\frac{\lambda N}{2}c^2
      +\frac N2\int_{\mathbb R^N}|u|^2\log|u|^2\,dx
      -\frac N2c^2\\
 &\quad+\mu\frac{N+\alpha}{2p}B_{\alpha,p}[u]
          +\frac{N+\beta}{2q}B_{\beta,q}[u]\\
={}&\frac{2-N}{2}\|\nabla u\|_2^2
      +\frac N2\bigl(\|\nabla u\|_2^2+\lambda c^2
                -\mu B_{\alpha,p}[u]-B_{\beta,q}[u]\bigr)\\
 &\quad-\frac{\lambda N}{2}c^2-\frac N2c^2
       +\mu\frac{N+\alpha}{2p}B_{\alpha,p}[u]
       +\frac{N+\beta}{2q}B_{\beta,q}[u]\\
={}&\|\nabla u\|_2^2-\frac N2c^2
       -\mu\left(\frac N2-\frac{N+\alpha}{2p}\right)B_{\alpha,p}[u]
       -\left(\frac N2-\frac{N+\beta}{2q}\right)B_{\beta,q}[u]\\
={}&\|\nabla u\|_2^2-\frac N2c^2
       -\mu\frac{Np-N-\alpha}{2p}B_{\alpha,p}[u]
       -\frac{Nq-N-\beta}{2q}B_{\beta,q}[u]\\
={}&\|\nabla u\|_2^2-\frac N2c^2
       -\mu\eta_pB_{\alpha,p}[u]-B_{\beta,q}[u]\\
={}&P(u),
\end{aligned}
\]
Here $(N-2)q=N+\beta$. Hence $(i)$ holds.

\medskip\noindent\textit{Part $(ii)$.}
We use an annular Pohozaev identity to show
$d(P|_{S_c})(u)\ne0$, then apply Lagrange multipliers
(cf. \cite[Lemma~4.2 and Theorem~4.1(iv)]{BartschSoave2017}). Let radial $w\in S_c$ solve
\[
 -a\Delta w+\lambda w
 =\delta w\log|w|^2
  +\mu b(I_\alpha*|w|^p)|w|^{p-2}w
  +d(I_\beta*|w|^q)|w|^{q-2}w,
\]
where $a,b,d,\lambda\in\mathbb R$ and $\delta\in\{0,1\}$,
including $a=0$. Then
\begin{equation}\label{eq:B2-generalized-pohozaev}
 a\|\nabla w\|_2^2-\frac{\delta N}{2}c^2
 -\mu b\eta_p B_{\alpha,p}[w]-dB_{\beta,q}[w]=0.
\end{equation}
Following the cutoff method of \cite[Proposition~3.1]{Moroz2013},
fix $0<\varepsilon<R/4$ and a radial cutoff $\chi$, equal to $1$
on $B_1$ and $0$ outside $B_2$. Put
\[
 V_{\varepsilon,R}(x)
 =x\chi(x/R)\bigl(1-\chi(x/\varepsilon)\bigr).
\]
On annuli, the radial lemma and the usual potential splitting ensure
local boundedness of $w$ and $I_\sigma*|w|^r$.
Thus the annular test is admissible, also for $a=0$, when the
Laplacian term is absent. Testing with $w$ gives
\begin{equation}\label{eq:B2-modified-test-u}
 a\|\nabla w\|_2^2+\lambda c^2
 =\delta\int_{\mathbb R^N}|w|^2\log|w|^2\,dx
  +\mu b B_{\alpha,p}[w]+d B_{\beta,q}[w].
\end{equation}
By direct calculation with $V_{\varepsilon,R}$, as
$\varepsilon\downarrow0$ and $R\uparrow\infty$, we obtain
\[
\begin{aligned}
 a\int_{\mathbb R^N}(-\Delta w)
             V_{\varepsilon,R}\cdot\nabla w\,dx
 &\longrightarrow\frac{(2-N)a}{2}\|\nabla w\|_2^2
 \qquad(a\ne0),\\
 \int_{\mathbb R^N}(\lambda w-\delta w\log|w|^2)
             V_{\varepsilon,R}\cdot\nabla w\,dx
 &\longrightarrow-\frac{\lambda N}{2}c^2
  +\frac{\delta N}{2}
       \left(\int_{\mathbb R^N}|w|^2\log|w|^2\,dx-c^2\right),\\
 \int_{\mathbb R^N}(I_\sigma*|w|^r)|w|^{r-2}w
             V_{\varepsilon,R}\cdot\nabla w\,dx
 &\longrightarrow-\frac{N+\sigma}{2r}B_{\sigma,r}[w],
 \qquad(\sigma,r)\in\{(\alpha,p),(\beta,q)\}.
\end{aligned}
\]
The errors vanish by integrability; the nonlocal limits follow
by HLS, symmetrization, and dominated convergence.
For $a=0$, the first term is simply omitted.
Together with \eqref{eq:B2-modified-test-u}, this yields
\[
\begin{aligned}
0={}&\frac{2-N}{2}a\|\nabla w\|_2^2
      -\frac{\lambda N}{2}c^2
      +\frac{\delta N}{2}
           \left(\int_{\mathbb R^N}|w|^2\log|w|^2\,dx-c^2\right)\\
 &\qquad+\mu b\frac{N+\alpha}{2p}B_{\alpha,p}[w]
                 +d\frac{N+\beta}{2q}B_{\beta,q}[w]\\
={}&\frac{2-N}{2}a\|\nabla w\|_2^2
     +\frac N2\left(
        a\|\nabla w\|_2^2+\lambda c^2
             -\mu bB_{\alpha,p}[w]-dB_{\beta,q}[w]\right)
       -\frac{\lambda N}{2}c^2-\frac{\delta N}{2}c^2\\
 &\qquad+\mu b\frac{N+\alpha}{2p}B_{\alpha,p}[w]
                   +d\frac{N+\beta}{2q}B_{\beta,q}[w]\\
={}&a\|\nabla w\|_2^2-\frac{\delta N}{2}c^2
       -\mu b\left(\frac N2-\frac{N+\alpha}{2p}\right)
                         B_{\alpha,p}[w]\\
 &\qquad-d\left(\frac N2-\frac{N+\beta}{2q}\right)
                         B_{\beta,q}[w]\\
={}&a\|\nabla w\|_2^2-\frac{\delta N}{2}c^2
       -\mu b\eta_pB_{\alpha,p}[w]-dB_{\beta,q}[w],
\end{aligned}
\]
We used $(N-2)q=N+\beta$, proving
\eqref{eq:B2-generalized-pohozaev}. If $d(P|_{S_c})(u)=0$, the multiplier rule gives
$\tau\in\mathbb R$ with
\[
\begin{aligned}
0&=P'(u)[v]+2\tau\int_{\mathbb R^N}uv\,dx\\
 &=2\int_{\mathbb R^N}\nabla u\cdot\nabla v\,dx
    -2\mu p\eta_p
      \int_{\mathbb R^N}(I_\alpha*|u|^p)|u|^{p-2}uv\,dx\\
 &\qquad-2q\int_{\mathbb R^N}
                  (I_\beta*|u|^q)|u|^{q-2}uv\,dx
     +2\tau\int_{\mathbb R^N}uv\,dx,
\end{aligned}
\]
for every $v\in H_r$. Radial averaging extends this identity
to smooth compactly supported tests. Hence
\[
 -\Delta u+\tau u
   =\mu p\eta_p(I_\alpha*|u|^p)|u|^{p-2}u
        +q(I_\beta*|u|^q)|u|^{q-2}u.
\]
Applying \eqref{eq:B2-generalized-pohozaev} with
$(a,\delta,b,d)=(1,0,p\eta_p,q)$ gives
\[
\begin{aligned}
0
 &=\|\nabla u\|_2^2
    -\mu p\eta_p^2B_{\alpha,p}[u]-qB_{\beta,q}[u]\\
 &=\frac12\left(2\|\nabla u\|_2^2
       -2\mu p\eta_p^2B_{\alpha,p}[u]
       -2qB_{\beta,q}[u]\right)\\
 &=\frac12Q(u),
\end{aligned}
\]
This contradicts $Q(u)\ne0$. Thus $d(P|_{S_c})(u)\ne0$,
and $\mathcal P_c$ is locally a $C^1$ hypersurface of $S_c$. If $d(E|_{\mathcal P_c})(u)=0$, the multiplier rule gives
$\lambda,\kappa\in\mathbb R$ such that
\[
 E'(u)[v]+\kappa P'(u)[v]
          +\lambda\int_{\mathbb R^N}uv\,dx=0
 \qquad(v\in H_r).
\]
Expanding and averaging radially yields
\[
\begin{aligned}
 -(1+2\kappa)\Delta u+\lambda u
 &=u\log|u|^2
   +\mu(1+2\kappa p\eta_p)
          (I_\alpha*|u|^p)|u|^{p-2}u\\
 &\qquad+(1+2\kappa q)
          (I_\beta*|u|^q)|u|^{q-2}u.
\end{aligned}
\]
Apply \eqref{eq:B2-generalized-pohozaev} with $a=1+2\kappa$, $\delta=1$, $b=1+2\kappa p\eta_p$, and $d=1+2\kappa q$, also for $\kappa=-\tfrac12$, when $a=0$. Then
\[
\begin{aligned}
0
 &=(1+2\kappa)\|\nabla u\|_2^2
     -\frac N2c^2
     -\mu(1+2\kappa p\eta_p)\eta_pB_{\alpha,p}[u]
     -(1+2\kappa q)B_{\beta,q}[u]\\
 &=\left(\|\nabla u\|_2^2-\frac N2c^2
       -\mu\eta_pB_{\alpha,p}[u]-B_{\beta,q}[u]\right)\\
 &\qquad+\kappa\left(2\|\nabla u\|_2^2
       -2\mu p\eta_p^2B_{\alpha,p}[u]
       -2qB_{\beta,q}[u]\right)\\
 &=P(u)+\kappa Q(u)\\
 &=\kappa Q(u).
\end{aligned}
\]
Since $Q(u)\ne0$, $\kappa=0$; hence
\[
 E'(u)+\lambda\Psi'(u)=0\quad\text{in }H_r^*,
 \qquad\text{hence}\qquad d(E|_{S_c})(u)=0.
\]
Thus $(ii)$ holds.

\end{proof}
}

\section{Existence of solutions}\label{sec:existence}

\begin{lemma}\label{lem 4.1}
Let $\mu > 0$ and let $u \in S_c$ satisfy $P(u) \le 0$. If $\|\nabla u\|_2^2 = \frac{K N c^2}{2}$ where $K > 1$ is a constant, then
\begin{align}\label{c-}
\frac{(K - 1)N c^2}{2} \le \mu \eta_p C_{N, \alpha, p}c^{2p}\left(\frac{K N}{2}\right)^{p\eta_p} + S_\beta^{-q}c^{2q}\left(\frac{K N}{2}\right)^q.
\end{align}
Therefore, if $u \in S_c$ satisfies $P(u) \le 0$ and
\begin{align}\label{c}
\frac{(K - 1)N c^2}{2} > \mu \eta_p C_{N, \alpha, p}c^{2p}\left(\frac{K N}{2}\right)^{p\eta_p} + S_\beta^{-q}c^{2q}\left(\frac{K N}{2}\right)^q.
\end{align}
then $\|\nabla u\|_2^2 \neq \frac{K N c^2}{2}$.
\end{lemma}

\begin{proof}
Let $u \in S_c$. We use the Pohozaev function \eqref{P}, apply the sharp Gagliardo-Nirenberg inequality \eqref{GN}, and we obtain
\begin{align}\label{4.1}
\|\nabla u\|_2^2 &\le \frac{N c^2}{2} + \mu \eta_p B_{\alpha, p}[u] + B_{\beta, q}[u]  \notag \\
& \le \frac{N c^2}{2} + |\mu| \eta_p C_{N, \alpha, p}c^{2p - 2p\eta_p}\|\nabla u\|_2^{2p\eta_p} + S_\beta^{-q}\|\nabla u\|_2^{2q}.
\end{align}
Substituting $\|\nabla u\|_2^2 = \frac{K N c^2}{2}$ in \eqref{4.1}, we obtain exactly \eqref{c-}. Therefore, if $u \in S_c$ satisfies $P(u) \le 0$ and \eqref{c} holds, then $\|\nabla u\|_2^2 \neq \frac{K N c^2}{2}$.
\end{proof}

\begin{remark}\label{rmk 2}
Condition \eqref{c} is equivalent to
\[
\frac{(K - 1)N}{2} > |\mu| \eta_p C_{N, \alpha, p}c^{2p - 2}\left(\frac{K N}{2}\right)^{p\eta_p} + S_\beta^{-q}c^{2q - 2}\left(\frac{K N}{2}\right)^q.
\]
Since $p>1$ and $q>1$, the exponents $2p - 2$ and $2q - 2$ are strictly positive. This condition is of the form $0 < c < c_K$, where $c_K > 0$ is the unique positive threshold determined by \eqref{c} and depends on the fixed parameters and $K>1$. This condition is called the small mass restriction and is used in various situations. In particular, in Lemma~\ref{lem 4.1}, it is used to exclude the possible solutions of \eqref{eq main} from the boundary of the open set
\[
W_{c, K}:= \left\{u \in S_c \mid \|\nabla u\|_2^2 < \frac{K N c^2}{2}\right\}.
\]
\end{remark}

In this subsection, we assume $\frac{N + \alpha}{N} < p < \frac{N + \alpha + 2}{N}$ and $\mu > 0$. As a result, $0 < 2p\eta_p < 2$. We proceed with the following lemmas:

\begin{lemma}\label{lem 4.3}
Let $c > 0$ satisfy \eqref{c} and let $K > \frac{q}{q - 1}$. Then $\mathcal{P}_c^0 = \varnothing$ provided that $0 < c < \min\{c_F, c_G\}$ where
\begin{align*}
c_F &=\sqrt{\frac{2}{N}}\left[\frac{S_\beta^q\left[(1 - p\eta_p)K + p\eta_p\right]}{(q - p\eta_p)K^q}\right]^{\frac{1}{2q - 2}},\\
c_G &=\left[\frac{N\left[(q - 1)K - q\right]}{2\mu\eta_p(q-p\eta_p)C_{N, \alpha, p}\left(\frac{KN}{2}\right)^{p\eta_p}}\right]^{\frac{1}{2p - 2}}.
\end{align*}
Moreover, there exists a constant $\delta > 0$ independent of $u$ such that $|Q(u)| \ge \delta > 0$, where the second Pohozaev functional $Q$ is defined in \eqref{Q}.
\end{lemma}

\begin{proof}
Let $u \in \mathcal{P}_c^0$. Then $P(u) = 0$ and $Q(u) = 0$. From \eqref{P} and \eqref{Q}, we obtain the following system:
\begin{align}\label{4.2}
\begin{cases}
\displaystyle{\|\nabla u\|_2^2 = \frac{N c^2}{2} + \mu \eta_p B_{\alpha, p}[u] + B_{\beta, q}[u]}, \\
\displaystyle{\|\nabla u\|_2^2 = \mu p\eta_p^2 B_{\alpha, p}[u] + qB_{\beta, q}[u]}.
\end{cases}
\end{align}
Solving \eqref{4.2}, we obtain
\begin{align}\label{4.3}
\begin{cases}
\displaystyle{(1 - p\eta_p)\|\nabla u\|_2^2 = -\frac{N p\eta_p c^2}{2} + (q - p\eta_p)B_{\beta, q}[u]},\\
\displaystyle{(q - 1)\|\nabla u\|_2^2 = \frac{N q c^2}{2} + \mu \eta_p(q - p\eta_p) B_{\alpha, p}[u]}.
\end{cases}
\end{align}
Applying the sharp inequalities \eqref{GS} and \eqref{S} on \eqref{4.3}, we get
\begin{align}\label{4.4}
\begin{cases}
\displaystyle{(1 - p\eta_p)\|\nabla u\|_2^2 \le -\frac{N p\eta_p c^2}{2} + (q - p\eta_p)S_\beta^{-q}\|\nabla u\|_2^{2q}},\\
\displaystyle{(q - 1)\|\nabla u\|_2^2 \le \frac{N q c^2}{2} + \mu \eta_p(q - p\eta_p)C_{N, \alpha, p}c^{2p - 2p\eta_p} \|\nabla u\|_2^{2p\eta_p}}.
\end{cases}
\end{align}
By the first inequality in \eqref{4.4}, we define the auxiliary function
\begin{align}\label{F}
F_c(t):= (1 - p\eta_p)t^2 + \frac{N p\eta_p c^2}{2} - (q - p\eta_p)S_\beta^{-q}t^{2q}, \qquad t > 0.
\end{align}
We aim to find conditions on $c, t > 0$ that guarantee $F_c(t) > 0$. Let $t_0:= \sqrt{\frac{K N c^2}{2}}$ satisfy $F_c(t_0) = 0$. Solving this equation on $c > 0$, we obtain
\[
c = c_F = \sqrt{\frac{2}{N}}\left[\frac{S_\beta^q\left[(1 - p\eta_p)K + p\eta_p\right]}{(q - p\eta_p)K^q}\right]^{\frac{1}{2q - 2}}.
\]
We observe that if $c = c_F$, then the value $t_0$ is the exact nonzero root of $F_c$. If $0 < c < c_F$, we have $F_c(t_0) > 0$. To see this, we can find
\[
\lim_{t \to 0^+} F_c(t) = \frac{N p\eta_p c^2}{2} > 0, \qquad  \lim_{t \to +\infty} F_c(t) = - \infty,
\]
and $F_c$ has exactly one critical point, which can be solved analytically and is a global maximum. Figure~\ref{fig:Fc} shows the graph of $F_c$.

\begin{figure}[H]
    \centering
    \includegraphics[width=0.75\textwidth]{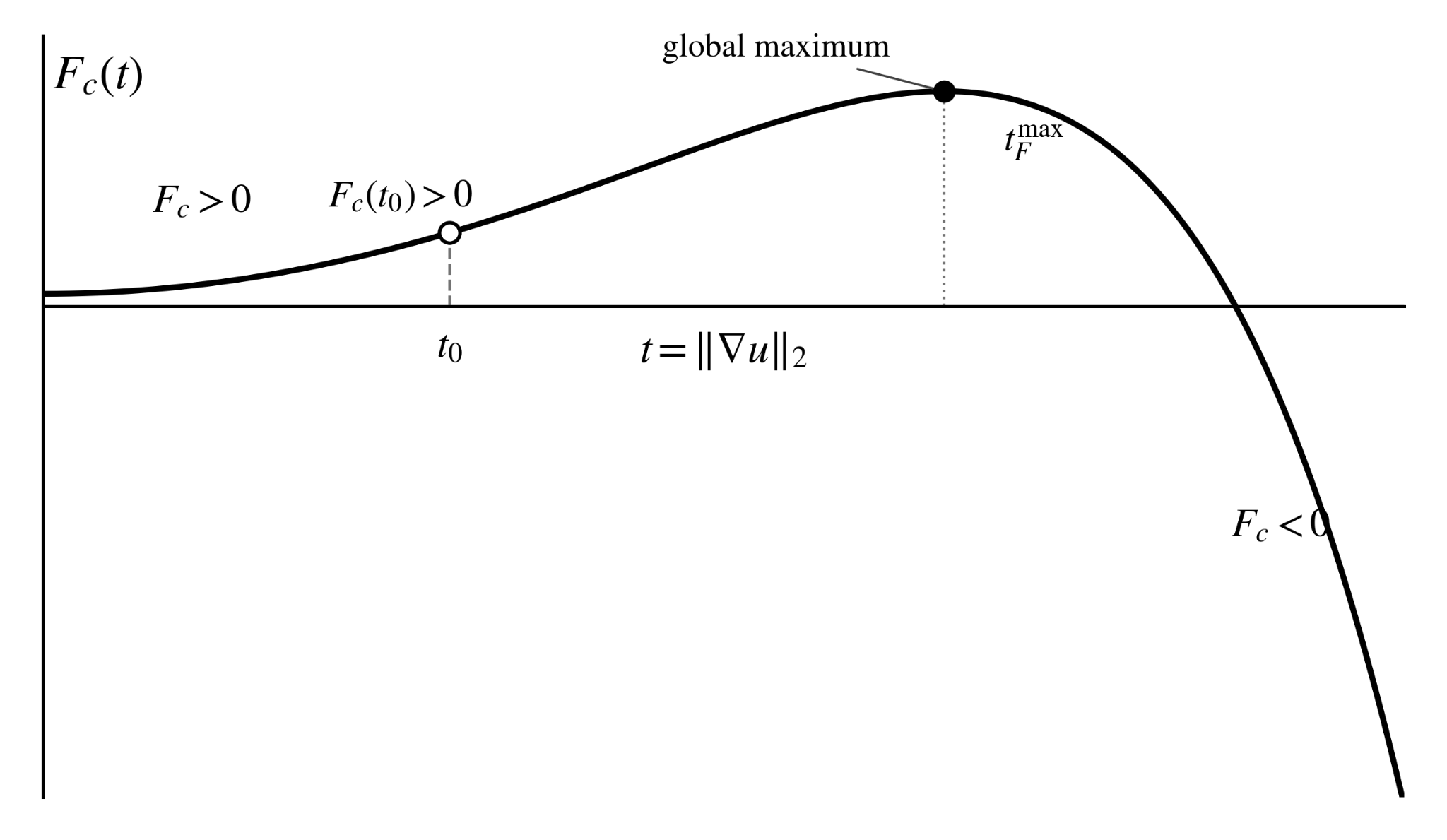}
    \caption{The geometric behavior of $F_c$ in the
    $L^2$-subcritical regime.}
    \label{fig:Fc}
\end{figure}

\noindent
{Let us indeed note that $F_c'(t)=2t[(1-p\eta_p)-q(q-p\eta_p)S_\beta^{-q}t^{2(q-1)}]$, which means that $F_c$ is strictly increasing until its (unique) critical point and strictly decreasing afterwards. Since  $F_c(0^+)>0$ and $F_c(t_0)>0$ for $0<c<c_F$, we deduce that $F_c(t)>0$ for all $0<t<t_0$.} Equivalently, $F_c(\|\nabla u\|_2)>0$ whenever $\|\nabla u\|_2<t_0$. Now, from the second inequality in \eqref{4.3}, we define the following auxiliary function
\begin{align}\label{G}
G_c(t):= (q - 1)t^2 - \frac{N q c^2}{2} - \mu \eta_p(q - p\eta_p)C_{N, \alpha, p}c^{2p - 2p\eta_p} t^{2p\eta_p}, \qquad t > 0.
\end{align}
We aim to find conditions on $c, t > 0$ that guarantee $G_c(t) > 0$. Recall $t_0= \sqrt{\frac{K N c^2}{2}}$. Let $t_0$ satisfy $G_c(t_0) = 0$. Solving this equation on $c > 0$, we obtain
\[
c = c_G = \left[\frac{N\left[(q - 1)K - q\right]}{2\mu\eta_p(q-p\eta_p)C_{N, \alpha, p}\left(\frac{KN}{2}\right)^{p\eta_p}}\right]^{\frac{1}{2p - 2}}.
\]
We observe that if $c = c_G$, then the value $t_0$ is the exact nonzero root of $G_c$. If $0 < c < c_G$, we have $G_c(t_0) > 0$. We also have
\[
\lim_{t \to 0^+} G_c(t) = - \frac{N q c^2}{2} < 0, \qquad  \lim_{t \to +\infty} G_c(t) = + \infty,
\]
and $G_c$ has exactly one critical point that can be solved analytically, which is a global minimum. Figure~\ref{fig:Gc} shows the graph of $G_c$.

\begin{figure}[H]
    \centering
    \includegraphics[width=0.75\textwidth]{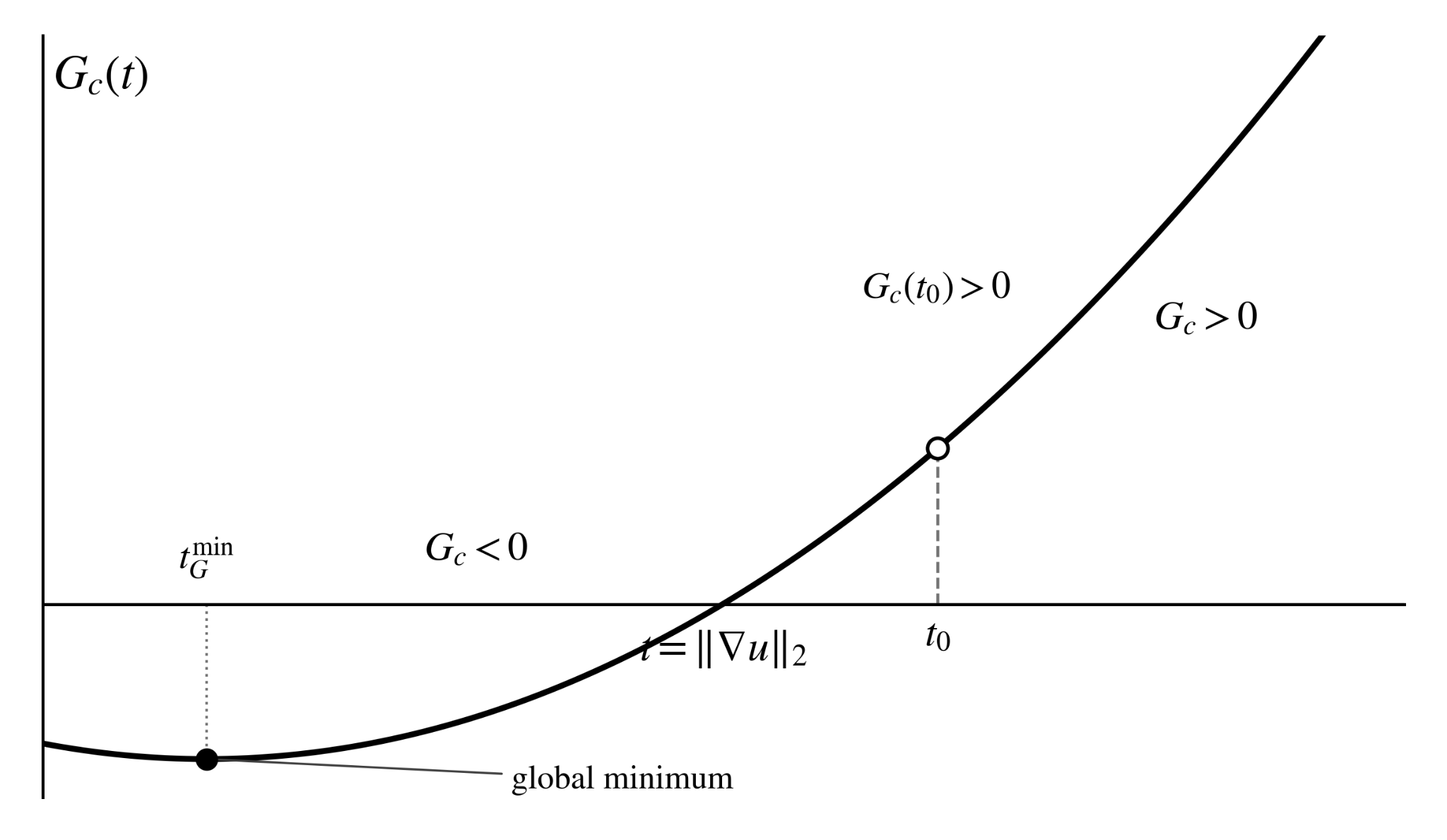}
    \caption{The geometric behavior of $G_c$ in the
    $L^2$-subcritical regime.}
    \label{fig:Gc}
\end{figure}

\noindent
{Again, we have
	\[
	G_c'(t)=2t\left[(q-1)-p\eta_p\,\mu\eta_p(q-p\eta_p)C_{N,\alpha,p}c^{2p-2p\eta_p}t^{2(p\eta_p-1)}\right].
	\]
	Since $0<p\eta_p<1$, the function $t\mapsto G_c'(t)/(2t)$ is strictly increasing on $(0,+\infty)$. Thus $G_c'$ vanishes exactly once, and $G_c$ admits a unique global minimum. Furthermore, $G_c(0^+)<0$ and $G_c(t_0)>0$ for $0<c<c_G$ imply that $t_0$ is located  to the right of this minimum and $G_c(t)>0$ for all $t>t_0$.} Equivalently, $G_c(\|\nabla u\|_2)>0$ whenever $\|\nabla u\|_2>t_0$. Therefore, if $0 < c < \min\{c_F, c_G\}$ we have $\mathcal{P}_c^0 = \varnothing$. Next, we prove the uniform bounds. Let $u \in \mathcal{P}_c$. By \eqref{P} and \eqref{Q}, after some calculation, we obtain
\begin{align}\label{4.5}
Q(u) \ge 2F_c(\|\nabla u\|_2), \qquad -Q(u) \ge 2G_c(\|\nabla u\|_2).
\end{align}
Let $t > 0$. We consider three cases. If $0 < t < t_0$, we obtain
\begin{align}\label{4.6}
F_c(t) &= (1 - p\eta_p)t^2 + \frac{N p\eta_p c^2}{2} - (q - p\eta_p)S_\beta^{-q}t^{2q} \notag \\
& \ge \frac{N p\eta_p c^2}{2} + t^2\left[(1 - p\eta_p) - (q - p\eta_p)S_\beta^{-q}t_0^{2q - 2}\right].
\end{align}
If the expression inside the bracket in the second inequality in \eqref{4.6} is positive, then $F_c(t) \ge \frac{N p\eta_p c^2}{2}$. Otherwise
\begin{align*}
F_c(t) &\ge \frac{N p\eta_p c^2}{2} + t^2\left[(1 - p\eta_p) - (q - p\eta_p)S_\beta^{-q}t_0^{2q - 2}\right] \\
& \ge \frac{N p\eta_p c^2}{2} + t_0^2\left[(1 - p\eta_p) - (q - p\eta_p)S_\beta^{-q}t_0^{2q - 2}\right] \\
& = F_c(t_0) > 0.
\end{align*}
To justify $F_c(t_0) > 0$, we refer to $0 < c < c_F$ and the following calculation
\begin{equation*}
F_c(t_0) = \frac{N}{2}\left[(1 - p\eta_p)K + p\eta_p\right]c^2\left[1-\left(\frac{c}{c_F}\right)^{2q - 2}\right] > 0.
\end{equation*}
By \eqref{4.5}, we get $Q(u) \ge \min\left\{\frac{N p\eta_p c^2}{2}, F_c(t_0)\right\} > 0$. If $t > t_0$, we then define
\begin{align}\label{4.7}
H_c(t) &:= \frac{G_c(t)}{t^2}= (q - 1) - \frac{N q c^2}{2t^2} - \mu \eta_p(q - p\eta_p)C_{N, \alpha, p}c^{2p - 2p\eta_p} t^{2p\eta_p - 2}.
\end{align}
By \eqref{4.7} and a direct calculation, we get
\[
H_c'(t) = \frac{N q c^2}{t^3} + 2\mu \eta_p(q - p\eta_p)(1 - p\eta_p)C_{N, \alpha, p}c^{2p - 2p\eta_p} t^{2p\eta_p - 3} > 0.
\]
Therefore $H_c$ is strictly increasing. Using this fact and noting that $t > t_0$, we obtain
\[
H_c(t) > H_c(t_0) \implies \frac{G_c(t)}{t^2} > \frac{G_c(t_0)}{t_0^2} > \frac{G_c(t_0)}{t^2} \implies G_c(t) > G_c(t_0) > 0.
\]
To justify $G_c(t_0) > 0$, we refer to $0 < c < c_G$ and the following calculation
\begin{equation*}
G_c(t_0)= \frac{N}{2}\left[(q - 1)K - q\right]c^2\left[1 - \left(\frac{c}{c_G}\right)^{2p - 2}\right],
\end{equation*}
and by \eqref{4.5}, we get $Q(u) \le -G_c(t_0) < 0$. The case $t = t_0$ is also impossible by Lemma~\ref{lem 4.1}. Combining these cases, we can find a $\delta > 0$ independent of $u$ (or $t$) such that $|Q(u)| \ge \delta > 0$. This completes the proof.
\end{proof}

\begin{remark}\label{rmk 4.4}
The small mass restrictions in Lemma~\ref{lem 4.1} and Lemma~\ref{lem 4.3} are indeed compatible. From Lemma~\ref{lem 4.1}, we have $0 < c < c_K$ and by Lemma~\ref{lem 4.3}, we have $0 < c < \min\{c_F, c_G\}$. Considering the intersection of intervals, we get $0 < \min\{c_F, c_G, c_K\}$ and all these numbers are nonzero. Therefore, the small mass conditions are compatible.
\end{remark}

\begin{lemma}\label{lem 4.5}
Let $0 < c < \min\{c_F, c_G, c_K\}$. For any $u \in S_c$, there exist numbers $-\infty < s^+ < s^- < +\infty$ such that $s^+ \star u \in \mathcal{P}_c^+$ and $s^- \star u \in \mathcal{P}_c^-$.
\end{lemma}

\begin{proof}
Recall from $s \star u(x) = e^{\frac{Ns}{2}}u(e^s x)$ and \eqref{fib} that
\begin{align}\label{4.8}
\Phi_u(s) &= \frac{e^{2s}}{2}\|\nabla u\|_2^2 + \frac{1}{2}\|u\|_2^2 - \frac{1}{2}\int_{\mathbb{R}^N}|u|^2\log|u|^2\,dx \notag \\
& \quad - \frac{N c^2 s}{2} - \frac{\mu e^{2p\eta_p s}}{2p}B_{\alpha, p}[u] - \frac{e^{2qs}}{2q}B_{\beta, q}[u].
\end{align}
Differentiating $\Phi_u$ from \eqref{4.8} in $s$, we obtain
\begin{align}\label{4.9}
\Phi_u'(s) = e^{2s}\|\nabla u\|_2^2 - \frac{N c^2}{2} - \mu \eta_p e^{2p\eta_p s}B_{\alpha, p}[u] - e^{2qs}B_{\beta, q}[u].
\end{align}
Dividing both sides of \eqref{4.9} by $e^{2p\eta_p s}$ and defining $\Psi_u(s):= \frac{\Phi_u'(s)}{e^{2p\eta_p s}}$, we obtain
\begin{align}\label{4.10}
\Psi_u(s) = e^{(2 - 2p\eta_p)s}\|\nabla u\|_2^2 - \frac{N c^2}{2e^{2p\eta_p s}} - \mu \eta_pB_{\alpha, p}[u] - e^{(2q - 2p\eta_p)s}B_{\beta, q}[u].
\end{align}
From \eqref{4.10}, we observe that
\[
\lim_{s \to -\infty} \Psi_u(s) = -\infty, \qquad \lim_{s \to +\infty} \Psi_u(s) = -\infty.
\]
Therefore, $\Psi_u$ has at least one maximum point. We show that this point is unique. Let $s^\ast$ be the critical point of $\Psi_u$. We can write
\begin{align}\label{4.11}
(2 - 2p\eta_p)\|\nabla u\|_2^2 + \frac{N p\eta_p c^2}{e^{2s^\ast}} = (2q - 2p\eta_p)e^{(2q - 2)s^\ast}B_{\beta, q}[u].
\end{align}
Since $0 < p\eta_p < 1 < q$, the left-hand side of \eqref{4.11} is strictly decreasing, while the right-hand side is strictly increasing, implying that the critical point is unique. By the structure, $\Psi_u$ has a global maximum. By Lemma~\ref{lem 4.1}, we have $P(s_0 \star u) > 0$ where $s_0 = \log\left(\frac{t_0}{\|\nabla u\|_2}\right)$. By \eqref{4.10}, we know that $\Psi_u(s) = e^{-2p\eta_p s}P(s \star u)$. Therefore, at $s_0$, we get $\Psi_u(s_0) > 0$. By the uniqueness of the critical point, $\Psi_u(s)$ has exactly two roots $-\infty < s^+ < s_0 < s^- < +\infty$ where $\Psi_u(s) > 0$ for $s^+ < s < s^-$ and $\Psi_u(s) < 0$ for $s < s^+$ or $s > s^-$. Since $\Phi_u'(s) = e^{2p\eta_p s}\Psi_u(s)$, the geometry is preserved for $\Phi_u'$ as well. By $s^+ < s_0 < s^-$, we get
\[
\|\nabla(s^+ \star u)\|_2 < t_0, \qquad \|\nabla(s^- \star u)\|_2 > t_0,
\]
and by Lemma~\ref{lem 4.3}, we obtain $Q(s^+ \star u) > 0$ and $Q(s^- \star u) < 0$, implying $s^+ \star u \in \mathcal{P}_c^+$ and $s^- \star u \in \mathcal{P}_c^-$.
\end{proof}

\begin{lemma}\label{lem 4.6} 
Let $\mu > 0$ and let $c > 0$ satisfy $0 < c < \min\{c_F, c_G, c_K\}$ and $K > \frac{q}{q - 1}$. Then, the value $m^+(c)$ is attained by a positive and radial function $u_c^+$ that is the solution to the problem \eqref{eq main} with the multiplier $\lambda_c^+ \in \mathbb{R}$, where
\[
m^+(c)= \inf_{u \in \mathcal{P}_c^+ \cap H_r} E(u).
\]
\end{lemma}

\begin{proof}
First, we show that $m^+(c) > -\infty$. Take a $u \in \mathcal{P}_c^+$ and consider $E(u) - \frac{1}{2q}P(u)$. Applying \eqref{GN} and \eqref{GS} and \eqref{S}, we obtain
\begin{align}\label{4.12}
E(u) & = E(u) - \frac{1}{2q}P(u) \notag \\
& = \left(\frac{1}{2} - \frac{1}{2q}\right)\|\nabla u\|_2^2 + \left(1 + \frac{N}{2q}\right)\frac{c^2}{2} \notag \\
& \quad - \frac{1}{2}\int_{\mathbb{R}^N} |u|^2\log{|u|^2}\,dx - \mu\left(\frac{1}{2p} - \frac{\eta_p}{2q}\right)B_{\alpha, p}[u] \notag \\
& \ge \left(\frac{1}{2} - \frac{1}{2q}\right)\|\nabla u\|_2^2 + \left(1 + \frac{N}{2q}\right)\frac{c^2}{2} \notag \\
& \quad - \frac{1}{2}\int_{\mathbb{R}^N} \mathscr{B}(|u|)\,dx - \mu\left(\frac{1}{2p} - \frac{\eta_p}{2q}\right)B_{\alpha, p}[u] \notag \\
& \ge \left(\frac{1}{2} - \frac{1}{2q}\right)\|\nabla u\|_2^2 + \left(1 + \frac{N}{2q}\right)\frac{c^2}{2} \notag \\
& - \frac{K_r}{2}\|u\|_r^r - \mu\left(\frac{1}{2p} - \frac{\eta_p}{2q}\right)C_{N, \alpha, p}c^{2p - 2p\eta_p}\|\nabla u\|_2^{2p\eta_p} \notag \\
& \ge \left(\frac{1}{2} - \frac{1}{2q}\right)\|\nabla u\|_2^2 + \left(1 + \frac{N}{2q}\right)\frac{c^2}{2} \notag \\
& - \frac{K_r C_{N, r}c^{r - r\gamma_r}}{2}\|\nabla u\|_2^{r\gamma_r} - \mu\left(\frac{1}{2p} - \frac{\eta_p}{2q}\right)C_{N, \alpha, p}c^{2p - 2p\eta_p}\|\nabla u\|_2^{2p\eta_p},
\end{align}
where $r > 2$ is taken so that $r\gamma_r < 2$. Also, by $0 < 2p\eta_p < 2$, we conclude that the lower bound in \eqref{4.12} remains finite for all $u \in \mathcal{P}_c^+$ (actually, it holds for all $u \in \mathcal{P}_c$). Therefore, $m^+(c) > -\infty$. To construct a Palais--Smale sequence for $m^+(c)$, we adapt the arguments in \cite[Lemma~5.6]{Shuai2025}. Since $\mathcal{P}_c^+$ is not a complete set, we apply Ekeland's variational principle directly on $S_c$. For any $u \in S_c$, by Lemma~\ref{lem 4.5}, we can find a unique $s^+ \in \mathbb{R}$ such that $s^+ \star u \in \mathcal{P}_c^+$. Following \cite{Shuai2025}, define the auxiliary map $J^+(u):= E(s^+ \star u)$ where $u \in S_c$. By the uniqueness of $s^+$, we obtain
\[
\inf_{u \in S_c}J^+(u) = \inf_{v \in \mathcal{P}_c^+}E(v) = m^+(c) > -\infty.
\]
By the Pohozaev identity \eqref{P}, for any $u \in \mathcal{P}_c$, we get
\[
\|\nabla u\|_2^2 = \frac{N c^2}{2} + \mu\eta_p B_{\alpha, p}[u] + B_{\beta, q}[u] \ge \frac{N c^2}{2} > 0,
\]
so the gradient norm is uniformly bounded from below. Now, let $(z_n)_{n \geq 1} \subset \mathcal{P}_c^+$ be a minimizing sequence such that
\[
m^+(c)\leq E(z_n)\leq m^+(c)+\frac{1}{n^2}.
\]
Since $z_n \in \mathcal{P}_c^+$, the estimate \eqref{4.12} yields $\sup_{n\geq 1}\|\nabla z_n\|_2<+\infty$. {The Pohozaev identity yields $\|\nabla z_n\|_2^2 \ge N c^2/2$. The sharp Gagliardo–Nirenberg and HLS inequalities give uniform boundedness of  the integral of $\mathscr B(|z_n|)$, $B_{\alpha,p}[z_n]$, and
	$B_{\beta,q}[z_n]$. Hence, we deduce that
	$\sup_n\int\mathscr A(|z_n|)<+\infty$ and so  $(z_n)$ is bounded in $H_r$.} Since $J^+(z_n) = E(z_n)$, applying Ekeland's variational principle to $J^+$ on the complete metric space $S_c$, we obtain a sequence $(w_n)_{n \geq 1} \subset S_c$ satisfying
\[
\begin{cases}
J^+(w_n)\to m^+(c),\\[2mm]
\|(J^+|_{S_c})'(w_n)\|\to 0,\\[2mm]
\|w_n-z_n\|_{H_r}\leq \dfrac{1}{n}.
\end{cases}
\]
Consequently, there exist constants $\delta, C >0$, independent of $n$, such that $\delta\leq\|\nabla w_n\|_2\leq C$, for all sufficiently large $n$. Now, set
\[
u_n=s_n^+\star w_n\in\mathcal{P}_c^+,
\qquad s_n^+=s^+(w_n).
\]
By the definition of $J^+$, we immediately have
\[
E(u_n)=J^+(w_n)\to m^+(c).
\]
Applying \eqref{4.12} to $u_n \in\mathcal{P}_c^+$, together with the Pohozaev identity \eqref{P}, we obtain $\delta\leq\|\nabla u_n\|_2\leq C$  after adjusting $\delta$ and $C$. Finally, the scaling identity implies
\[
e^{2s_n^+}
=\frac{\|\nabla u_n\|_2^2}
       {\|\nabla w_n\|_2^2}.
\]
Hence, $\sup_{n\geq 1}|s_n^+|<+\infty$. Also, by using the uniform boundedness of the norms in \eqref{orlicz-recovery-identity}, we conclude that
\[
\left|\int_{\mathbb{R}^N} |u_n|^2\log |u_n|^2\,dx\right| \le \int_{\mathbb{R}^N}\mathscr{A}(|u_n|)\,dx + \int_{\mathbb{R}^N}\mathscr{B}(|u_n|)\,dx, \le C.
\]
Since $Q(s^+(w)\star w)>0$, the implicit function theorem gives
$s^+\in C^1(S_c,\mathbb R)$. For $h\in T_{w_n}S_c$, the chain rule and $P(u_n)=0$ yield
\[
(J^+|_{S_c})'(w_n)[h]
=E'(u_n)[s_n^+\star h]
+P(u_n)(s^+)'(w_n)[h]
=E'(u_n)[s_n^+\star h].
\]
For $v\in T_{u_n}S_c$, put $h=(-s_n^+)\star v\in T_{w_n}S_c$. Since
$\sup_n|s_n^+|<\infty$, the mass-preserving dilation and its inverse are
uniformly bounded on $H_r$; hence
\[
\begin{aligned}
|E'(u_n)[v]|
&=|(J^+|_{S_c})'(w_n)[(-s_n^+)\star v]|
\le C\|(J^+|_{S_c})'(w_n)\|\,\|v\|_{H_r}.
\end{aligned}
\]
Thus $(u_n)\subset\mathcal P_c^+$ satisfies
\[
E(u_n)\to m^+(c),\qquad
\|(E|_{S_c})'(u_n)\|\to0.
\]
Since $(u_n)_{n \ge 1} \subset \mathcal{P}_c^+$ is a PS$_c$ sequence, we can apply the Lagrange multiplier rule and find a sequence $(\lambda_n)_{n \ge 1} \in \mathbb{R}$ such that
\[
E'(u_n) + \lambda_n u_n \to 0, \quad \text{in } {H_r^\ast,}
\]
as $n \to +\infty$. Therefore
\begin{align}\label{4.14}
-\Delta u_n + \lambda_n u_n &= u_n\log|u_n|^2 + \mu(I_\alpha \ast |u_n|^p)|u_n|^{p - 2}u_n \notag \\
& \quad + (I_\beta \ast |u_n|^q)|u_n|^{q - 2}u_n + o_n(1), \quad \text{in } H_r^\ast.
\end{align}
Testing \eqref{4.14} with $u_n$, we arrive at the identity
\[
\|\nabla u_n\|_2^2 + \lambda_n c^2 = \int_{\mathbb{R}^N} |u_n|^2\log |u_n|^2\,dx + \mu B_{\alpha, p}[u_n] + B_{\beta, q}[u_n] + o_n(1).
\]
Using the boundedness of all norms, we conclude that $(\lambda_n)_{n \ge 1}$ is uniformly bounded; there exists a $\lambda \in \mathbb{R}$ such that, up to a subsequence, $\lambda_n \to \lambda$ as $n \to +\infty$.
{Moreover, the boundedness of $(u_n)$ in the reflexive Banach space $H_r$ implies that there exists a $u \in H_r$ such that
\[
\begin{cases}
	u_n \rightharpoonup u, \qquad \text{weakly in } H_r, \\
	u_n \to u, \qquad \text{a.e. in } \mathbb{R}^N, \\
	u_n \to u, \qquad \text{strongly in } L^r(\mathbb{R}^N), \qquad \text{for } 2 \le r < \frac{2N}{N - 2},
\end{cases}
\]
where the strong convergence in $L^2(\mathbb R^N)$ follows by Lemma~\ref{lem 2.1}, while the strong convergence in
$L^s(\mathbb R^N)$ for all $2\le s<2^*$ follows by
interpolation with the boundedness of $(u_n)$ in
$L^{2^*}(\mathbb R^N)$.	Consequently,	the sequence $(|u_n|^s)$ converges strongly in $L^1(\mathbb R^N)$.	It follows that
$(\mathscr B(|u_n|))$ has uniformly absolutely continuous
integrals. Hence, by Vitali's theorem,
\[
\int_{\mathbb R^N}\mathscr B(|u_n|)\,dx
\longrightarrow
\int_{\mathbb R^N}\mathscr B(|u|)\,dx.
\]
By Lemma~\ref{lem 2.5}, the boundedness and weak convergence of $(u_n)$ in $H_r$,
as well as strong convergence in
$L^2(\mathbb R^N)$ imply
$
B_{\alpha,p}[u_n]\longrightarrow B_{\alpha,p}[u].
$	
For the subcritical Choquard term, by Lemma~\ref{lem 2.5} we conclude convergence against all test functions in $H_r$. For the critical Choquard term, the strong convergence of $(u_n)$ and the HLS inequality give convergence against every smooth test function in $C_c^\infty(\mathbb R^N)$.\\
Passing to the weak limits in \eqref{4.14}, we conclude that $u \in H_r$ solves it. Now, the principle of symmetric criticality implies that every radial normalized critical point is a weak solution in the full test space $H^\ast$. Then, }
\begin{align}\label{4.15}
-\Delta u + \lambda u &= u\log|u|^2 + \mu(I_\alpha \ast |u|^p)|u|^{p - 2}u \notag \\
& \quad + (I_\beta \ast |u|^q)|u|^{q - 2}u, \quad \text{in } H^\ast.
\end{align}
Since $u_n\in\mathcal P_c^+$, we have $P(u_n)=0$.
By \eqref{4.15} and Lemma~\ref{lem 2.4}(i), we get $P(u)=0$.
It remains to prove the strong convergence. We show that $u \in \mathcal{P}_c^+$. By the lower semi-continuity of the gradient norm, we get
\[
\|\nabla u\|_2^2 \le \liminf_{n \to +\infty} \|\nabla u_n\|_2^2 \le \frac{K N c^2}{2} \implies \|\nabla u\|_2^2 \le \frac{K N c^2}{2},
\]
and by Lemma~\ref{lem 4.1}, $\|\nabla u\|_2^2 \neq \frac{K N c^2}{2}$. Therefore $\|\nabla u\|_2^2 < \frac{K N c^2}{2}$. Hence by Lemma~\ref{lem 4.3}, $u \in \mathcal{P}_c^+$. Applying the convergence results, using the lower semi-continuity properties
\[
\|\nabla u\|_2^2 \le \liminf_{n \to +\infty}\|\nabla u_n\|_2^2, \qquad \int_{\mathbb{R}^N} \mathscr{A}(|u|)\,dx \le \liminf_{n \to +\infty} \int_{\mathbb{R}^N} \mathscr{A}(|u_n|)\,dx,
\]
and calculating similar to \eqref{4.12}, we obtain
\begin{align}\label{4.16}
m^+(c) &\le E(u) = E(u) - \frac{1}{2q}P(u) \notag \\
& = \left(\frac{1}{2} - \frac{1}{2q}\right)\|\nabla u\|_2^2 + \left(1 + \frac{N}{2q}\right)\frac{\|u\|_2^2}{2} \notag \\
& \quad + \frac{1}{2}\int_{\mathbb{R}^N} \mathscr{A}(|u|)\,dx - \frac{1}{2}\int_{\mathbb{R}^N} \mathscr{B}(|u|)\,dx - \mu\left(\frac{1}{2p} - \frac{\eta_p}{2q}\right)B_{\alpha, p}[u] \notag \\
& \le \liminf_{n \to +\infty}\Bigg[\left(\frac{1}{2} - \frac{1}{2q}\right)\|\nabla u_n\|_2^2 + \left(1 + \frac{N}{2q}\right)\frac{\|u_n\|_2^2}{2} \notag \\
& \quad + \frac{1}{2}\int_{\mathbb{R}^N} \mathscr{A}(|u_n|)\,dx - \frac{1}{2}\int_{\mathbb{R}^N} \mathscr{B}(|u_n|)\,dx - \mu\left(\frac{1}{2p} - \frac{\eta_p}{2q}\right)B_{\alpha, p}[u_n]\Bigg] \notag \\
& = \liminf_{n \to +\infty}\left[E(u_n) - \frac{1}{2q}P(u_n)\right] \notag \\
& = \liminf_{n \to +\infty} E(u_n) \notag \\
& = m^+(c).
\end{align}
From \eqref{4.16} we conclude that $m^+(c) = E(u)$. Thus, $m^+(c)$ is attained by $u$. Continuing with $P(u_n) = 0$, $P(u) = 0$, the multiplier identities, and $m^+(c) = E(u)$, we get
\begin{align}\label{4.17}
\begin{cases}
\|\nabla u_n\|_2^2 = \frac{N}{2}\|u_n\|_2^2 + \mu\eta_p B_{\alpha, p}[u_n] + B_{\beta, q}[u_n], \\
\|\nabla u\|_2^2 = \frac{N}{2}\|u\|_2^2 + \mu\eta_p B_{\alpha, p}[u] + B_{\beta, q}[u], \\
\|\nabla u_n\|_2^2 + \lambda_n \|u_n\|_2^2 = \displaystyle\int_{\mathbb{R}^N} |u_n|^2\log |u_n|^2\,dx + \mu B_{\alpha, p}[u_n] + B_{\beta, q}[u_n] + o_n(1), \\
\|\nabla u\|_2^2 + \lambda \|u\|_2^2 = \displaystyle\int_{\mathbb{R}^N} |u|^2\log |u|^2\,dx + \mu B_{\alpha, p}[u] + B_{\beta, q}[u], \\
m^+(c) = \left(\frac{1}{2} - \frac{1}{2q}\right)\|\nabla u\|_2^2 + \left(1 + \frac{N}{2q}\right)\frac{\|u\|_2^2}{2} \\
 \quad - \frac{1}{2}\displaystyle\int_{\mathbb{R}^N}|u|^2\log |u|^2\,dx - \mu\left(\frac{1}{2p} - \frac{\eta_p}{2q}\right)B_{\alpha, p}[u], \\
m^+(c) + o_n(1) = \left(\frac{1}{2} - \frac{1}{2q}\right)\|\nabla u_n\|_2^2 + \left(1 + \frac{N}{2q}\right)\frac{\|u_n\|_2^2}{2} \\
 \quad - \frac{1}{2}\displaystyle\int_{\mathbb{R}^N}|u_n|^2\log |u_n|^2\,dx - \mu\left(\frac{1}{2p} - \frac{\eta_p}{2q}\right)B_{\alpha, p}[u_n].
\end{cases}
\end{align}
and passing to the limit in \eqref{4.17}, we obtain
\begin{align}\label{4.18}
\begin{cases}
\lim_{n \to +\infty} \left[\|\nabla u_n\|_2^2 - B_{\beta, q}[u_n]\right] = \|\nabla u\|_2^2 - B_{\beta, q}[u], \\
\lim_{n \to +\infty} \left[\|\nabla u_n\|_2^2 - \displaystyle\int_{\mathbb{R}^N} |u_n|^2\log |u_n|^2\,dx - B_{\beta, q}[u_n]\right] \\
= \|\nabla u\|_2^2 - \displaystyle\int_{\mathbb{R}^N} |u|^2\log |u|^2\,dx - B_{\beta, q}[u], \\
\lim_{n \to +\infty} \Bigg[\left(\frac{1}{2} - \frac{1}{2q}\right)\|\nabla u_n\|_2^2 - \frac{1}{2}\displaystyle\int_{\mathbb{R}^N}|u_n|^2\log |u_n|^2\,dx \\
= \left(\frac{1}{2} - \frac{1}{2q}\right)\|\nabla u\|_2^2 - \frac{1}{2}\displaystyle\int_{\mathbb{R}^N}|u|^2\log |u|^2\,dx\Bigg].
\end{cases}
\end{align}
Solving \eqref{4.18}, we get
\[
\lim_{n \to +\infty} \int_{\mathbb{R}^N}|u_n|^2\log |u_n|^2\,dx = \int_{\mathbb{R}^N}|u|^2\log |u|^2\,dx.
\]
and
\[
\lim_{n \to +\infty} \|\nabla u_n\|_2^2 = \|\nabla u\|_2^2, \qquad \lim_{n \to +\infty} B_{\beta, q}[u_n] = B_{\beta, q}[u],
\]
By the definition $|u_n|^2\log |u_n|^2 = \mathscr{B}(|u_n|) - \mathscr{A}(|u_n|)$, we get $\mathscr{A}(|u_n|) \to \mathscr{A}(|u|)$ strongly in $L^1(\mathbb{R}^N)$. Notice that the strong convergence and Lemma~\ref{lem 4.3} yield
\[
Q(u) = \lim_{n \to +\infty} Q(u_n) \ge \delta' > 0,
\]
where $\delta' > 0$ is a constant independent of $n$. Therefore, $Q(u)$ is uniformly positive. Finally, we prove that $u$ is radial and strictly positive. The first statement is straightforward by working on $H_r$. For the second, we observe that \(E(|u|)=E(u)=m^+(c)\). Therefore, we can assume that $|u|$ is also a minimizer. Hence $u \ge 0$. To prove strict positivity, rewrite \eqref{eq main} as
\[
-\Delta u + \lambda u - u\log|u|^2 = \mu(I_\alpha \ast |u|^p)|u|^{p - 2}u+(I_\beta \ast |u|^q)|u|^{q - 2}u \ge 0.
\]
Therefore $-\Delta u + f(u) \ge 0$, where $f(t) = \lambda t - t\log |t|^2$ for $t>0$ and  $f(0) = 0$. We check Vázquez's conditions. For $t > 0$ small enough, we have
\[
f'(t) = \lambda - \log |t|^2 - 2 > 0.
\]
Also
\[
\int_0^\eta \frac{ds}{\sqrt{s f(s)}} = \int_0^\eta \frac{ds}{s\sqrt{ \lambda - \log |s|^2}} = +\infty,
\]
for $\eta > 0$ small enough. Since $\|u\|_2 = c > 0$, the solution $u$ cannot be identically zero, and Lemma~\ref{lem 2.6} guarantees the required regularity for applying Vázquez's strong maximum principle \cite{Vazquez1984}, and we get $u(x) > 0$ for all $x \in \mathbb{R}^N$. We denote this solution by $u_c^+$ and its Lagrange
multiplier by $\lambda_c^+$. This completes the proof.
\end{proof}

\noindent
Now, we construct the second solution for problem \eqref{eq main}.

\begin{lemma}\label{lem 4.7}
Assume that $\mu>0$, $K>\frac{q}{q-1}$, and
$0<c<\min\{c_F,c_G,c_K\}$. Then
\[
m^+(c)=\inf_{v\in W_{c,K}}E(v)
<\inf_{v\in\partial W_{c,K}}E(v)
\le m^-(c)<+\infty,
\]
where
\[
W_{c,K}:=\left\{v\in S_c:\ \|\nabla v\|_2^2
<\frac{KNc^2}{2}\right\}.
\]
In particular, $m^-(c)>m^+(c)$.
\end{lemma}

\begin{proof}
By Lemma~\ref{lem 4.1} and the small-mass condition, there is a
constant $\delta_0>0$ such that
\begin{equation}\label{4.19}
P(v)\ge\delta_0>0\qquad
\text{whenever }v\in\partial W_{c,K}.
\end{equation}
Indeed, the sharp Gagliardo--Nirenberg and HLS inequalities give
\[
\begin{aligned}
P(v)
&\ge\|\nabla v\|_2^2-\frac{Nc^2}{2}
-\mu\eta_p C_{N,\alpha,p}c^{2p-2p\eta_p}
\|\nabla v\|_2^{2p\eta_p}
-S_\beta^{-q}\|\nabla v\|_2^{2q}\\
&=\frac{(K-1)Nc^2}{2}
-\mu\eta_p C_{N,\alpha,p}c^{2p}
\left(\frac{KN}{2}\right)^{p\eta_p}
-S_\beta^{-q}c^{2q}\left(\frac{KN}{2}\right)^q>0
\end{aligned}
\]
on $\partial W_{c,K}$. Moreover, for such $v$ and $s\le0$,
\begin{equation}\label{4.20}
\begin{aligned}
P(s\star v)\ge{}&e^{2s}\frac{KNc^2}{2}-\frac{Nc^2}{2}\\
&-\mu\eta_p C_{N,\alpha,p}c^{2p}
\left(\frac{KN}{2}\right)^{p\eta_p}e^{2p\eta_p s}\\
&-S_\beta^{-q}c^{2q}\left(\frac{KN}{2}\right)^q e^{2qs}.
\end{aligned}
\end{equation}
Since the right-hand side of \eqref{4.20} is strictly positive
at $s=0$, there exist $\varepsilon_0,\delta_1>0$, independent of
$v\in\partial W_{c,K}$, such that
\begin{equation}\label{4.21}
P(s\star v)\ge\delta_1\qquad
\text{for all }s\in[-\varepsilon_0,0].
\end{equation}
By Lemma~\ref{lem 4.5}, every $v\in S_c$ has unique fibering
parameters $s^+(v)<s^-(v)$, and
\[
s^+(v)<\log\frac{\sqrt{KNc^2/2}}{\|\nabla v\|_2}<s^-(v).
\]
For $v\in\partial W_{c,K}$, these inequalities imply
$s^+(v)<0<s^-(v)$. In view of \eqref{4.21},
$s^+(v)\le-\varepsilon_0$, and hence
\begin{equation}\label{4.22}
\begin{aligned}
E(v)
&=E(s^+(v)\star v)+\int_{s^+(v)}^0P(s\star v)\,ds\\
&\ge m^+(c)+\int_{-\varepsilon_0}^0P(s\star v)\,ds\\
&\ge m^+(c)+\varepsilon_0\delta_1.
\end{aligned}
\end{equation}
Consequently,
$\inf_{v\in\partial W_{c,K}}E(v)>m^+(c)$. If $v\in W_{c,K}$, then $s^-(v)>0$. The fibering geometry in
Lemma~\ref{lem 4.5} shows that
\begin{equation}\label{4.23}
E(v)\ge E(s^+(v)\star v)\ge m^+(c),
\end{equation}
regardless of whether $s^+(v)$ is positive or nonpositive.
Since $u_c^+\in W_{c,K}$ and $E(u_c^+)=m^+(c)$, it follows that
\[
m^+(c)=\inf_{v\in W_{c,K}}E(v).
\]
Finally, let $v\in\mathcal P_c^-$. By Lemma~\ref{lem 4.5},
$\|\nabla v\|_2^2>KNc^2/2$ and
\[
s^+(v)<\log\frac{\sqrt{KNc^2/2}}{\|\nabla v\|_2}
<s^-(v)=0.
\]
The function $s\mapsto E(s\star v)$ is strictly increasing
between $s^+(v)$ and $s^-(v)$. Therefore
\begin{equation}\label{4.24}
E(v)>E\left(\left(\log\frac{\sqrt{KNc^2/2}}
{\|\nabla v\|_2}\right)\star v\right)
\ge\inf_{w\in\partial W_{c,K}}E(w).
\end{equation}
Taking the infimum over $\mathcal P_c^-$ proves
$m^-(c)\ge\inf_{\partial W_{c,K}}E>m^+(c)$.
The finiteness of $m^-(c)$ follows from the nonemptiness of
$\mathcal P_c^-$ established in Lemma~\ref{lem 4.5}.
\end{proof}

\begin{lemma}\label{lem 4.8}
Assume that $K>\frac{q}{q-1}$ and
$0<c<\min\{c_F,c_G,c_K\}$. Then there exists a sequence
$\{u_n\}\subset\mathcal P_c^-$ such that
\[
\begin{cases}
E(u_n)\to m^-(c),\\[1mm]
\big\|(E|_{S_c})'(u_n)\big\|_{(T_{u_n}S_c)^*}
\to0,\\[1mm]
P(u_n)=0,\\
 Q(u_n)\le-\delta<0,
\end{cases}
\]
where $\delta>0$ is independent of $n$.
Moreover, $\{u_n\}$ is bounded in $H_r$.
\end{lemma}

\begin{proof}
Choose $2<r<2^*$ with $r\gamma_r<2$. For $v\in\mathcal P_c$, \eqref{4.12} and Young's inequality give
\[
 E(v)=E(v)-\frac1{2q}P(v)
 \ge \frac{q-1}{4q}\|\nabla v\|_2^2-C_c.
\]
Thus $m^-(c)>-\infty$ by Lemma~\ref{lem 4.5}. For $v\in\mathcal P_c$ with bounded energy, the preceding bound gives bounded $\|\nabla v\|_2$. By \eqref{GN}, \eqref{GS}, \eqref{ineq}, and \eqref{orlicz-recovery-identity},
\[
 0\le\int_{\mathbb R^N}\mathscr A(|v|)\,dx
 \le 2E(v)+C_c\bigl(1+\|\nabla v\|_2^{r\gamma_r}
       +\|\nabla v\|_2^{2p\eta_p}
       +\|\nabla v\|_2^{2q}\bigr).
\]
Property~$(i)$ of $G$ then yields $\|v\|_{H_r}\le C$ on every bounded-energy subset of $\mathcal P_c$.

For $v\in S_c$, let $s^-(v)$ be given by Lemma~\ref{lem 4.5}. Set $J^-(v)=E(s^-(v)\star v)$; then $\inf_{S_c}J^-=m^-(c)$.
The identity $P(s^-(v)\star v)=0$ and Lemma~\ref{lem 4.3} give
\[
 \partial_s P(s\star v)\big|_{s=s^-(v)}
 =Q(s^-(v)\star v)\le-\delta<0.
\]
The implicit-function theorem and \eqref{fib} imply $s^-\in C^1(S_c)$ and $J^-\in C^1(S_c)$. For $h\in T_vS_c$,
\[
 (J^-)'(v)[h]
 =E'(s^-(v)\star v)[s^-(v)\star h].
\]
Indeed, the $s^-$-derivative vanishes since $P(s^-(v)\star v)=0$.

Take $z_n\in\mathcal P_c^-$ with
$ m^-(c)\le E(z_n)\le m^-(c)+n^{-2}$.
Then $J^-(z_n)=E(z_n)$ and $\|z_n\|_{H_r}\le C$.
Ekeland's principle on the complete metric space $S_c$ gives $w_n\in S_c$ such that
\[
\begin{cases}
 J^-(w_n)\to m^-(c),\\
 \|(J^-)'(w_n)\|_{(T_{w_n}S_c)^*}\le \frac{1}{n},\\
 \|w_n-z_n\|_{H_r}\le \frac{1}{n}.
\end{cases}
\]
As $z_n\in\mathcal P_c^-$, Lemma~\ref{lem 4.5} gives $\|\nabla z_n\|_2>t_0$. Hence $0<a\le\|\nabla w_n\|_2\le b$ for $n\gg1$. Put $s_n=s^-(w_n)$ and $u_n=s_n\star w_n\in\mathcal P_c^-$. Then
$E(u_n)=J^-(w_n)\to m^-(c)$, and the first estimate gives $\|u_n\|_{H_r}\le C$. Moreover,
\[
 t_0<\|\nabla u_n\|_2\le C,\qquad
 e^{2s_n}=\frac{\|\nabla u_n\|_2^2}
                   {\|\nabla w_n\|_2^2},\qquad
 |s_n|\le C.
\]
The $L^2$-preserving dilation and its inverse are uniformly bounded on $H_r$ for bounded $|s_n|$ (by the $\Delta_2$ property of $\mathscr A$). For $v\in T_{u_n}S_c$, let $h=(-s_n)\star v\in T_{w_n}S_c$. By the derivative identity above,
\[
 |(E|_{S_c})'(u_n)[v]|
 =|(J^-)'(w_n)[h]|
 \le\frac{C}{n}\|v\|_{H_r}.
\]
Consequently,
\[
\begin{cases}
 E(u_n)\to m^-(c),\\
 \|(E|_{S_c})'(u_n)\|_{(T_{u_n}S_c)^*}\to0,\\
 P(u_n)=0,\\
 Q(u_n)\le-\delta<0.
\end{cases}
\]
The proof is complete.
\end{proof}

\noindent
Let $N \ge 3$, $0 < \alpha, \beta < N$, and $q = \frac{N + \beta}{N - 2}$. Let $\kappa_{N, \beta} > 0$ be chosen such that the function
\[
U_1(x) = \kappa_{N, \beta}(1 + |x|^2)^{-\frac{N - 2}{2}},
\]
is the extremal function for the HLS inequality and
\[
\|\nabla U_1\|_2^2 = B_{\beta, q}[U_1] = S_\beta^{\frac{N + \beta}{\beta + 2}}.
\]
For $0 < \varepsilon < 1$ and $N = 3, 4$ we define the cutoff bubble following \cite{ChenJinRadulescuZhou2026} as follows:
\begin{align}\label{4.25}
U_\varepsilon(x)=
\kappa_{N,\beta}
\begin{cases}
\displaystyle
\left(\frac{\varepsilon}{\varepsilon^2+|x|^2}\right)^{\frac{N-2}{2}},
& |x|<1,\\[3mm]
\displaystyle
\left(\frac{\varepsilon}{1+\varepsilon^2}\right)^{\frac{N-2}{2}}
(2-|x|),
&1\le |x|<2,\\[3mm]
0,&|x|\ge2,
\end{cases}
\end{align}
and for $N \ge 5$, we use the full bubble
\begin{align}\label{4.26}
U_\varepsilon(x)=\kappa_{N,\beta}\left(\frac{\varepsilon}{\varepsilon^2+|x|^2}\right)^{\frac{N-2}{2}}.
\end{align}

\begin{lemma}\label{lem 4.9}
Assume that $K > \frac{q}{q - 1}$ and $0 < c < \min\{c_F, c_G, c_K\}$. Then $U_\varepsilon \in H_r$. If $N = 3$, then
\begin{align}\label{4.27}
\begin{aligned}
\|\nabla U_\varepsilon\|_2^2
&=
S_\beta^{\frac{N+\beta}{\beta+2}}
+O(\varepsilon),\\
B_{\beta,q}[U_\varepsilon]
&=
S_\beta^{\frac{N+\beta}{\beta+2}}
+O\left(\varepsilon^{\frac{N+\beta}{2}}\right),\\
\|U_\varepsilon\|_2^2
&=
\frac{23}{15}\kappa_{N,\beta}^2
|\mathbb S^2|\varepsilon
+O(\varepsilon^2),\\
\int_{\mathbb R^3}U_\varepsilon^2
\log U_\varepsilon^2\,dx
&=
\frac{23}{15}\kappa_{N,\beta}^2
|\mathbb S^2|\varepsilon\log\varepsilon
+O(\varepsilon).
\end{aligned}
\end{align}
If $N = 4$, then
\begin{align}\label{4.28}
\begin{aligned}
\|\nabla U_\varepsilon\|_2^2
&=
S_\beta^{\frac{N+\beta}{\beta+2}}
+O(\varepsilon^2),\\
B_{\beta,q}[U_\varepsilon]
&=
S_\beta^{\frac{N+\beta}{\beta+2}}
+O\left(\varepsilon^{\frac{N+\beta}{2}}\right),\\
\|U_\varepsilon\|_2^2
&=
\kappa_{N,\beta}^2|\mathbb S^3|
\varepsilon^2|\log\varepsilon|
+O(\varepsilon^2),\\
\int_{\mathbb R^4}U_\varepsilon^2
\log U_\varepsilon^2\,dx
&=
\kappa_{N,\beta}^2|\mathbb S^3|
\left(\log\kappa_{N,\beta}^2-\frac{12}{5}\right)
\varepsilon^2|\log\varepsilon|
+O(\varepsilon^2).
\end{aligned}
\end{align}
If $N \ge 5$, then
\begin{align}\label{4.29}
\begin{aligned}
\|\nabla U_\varepsilon\|_2^2
&=S_\beta^{\frac{N+\beta}{\beta+2}},\\
B_{\beta,q}[U_\varepsilon]
&=S_\beta^{\frac{N+\beta}{\beta+2}},\\
\|U_\varepsilon\|_2^2
&=\varepsilon^2\|U_1\|_2^2,\\
\int_{\mathbb R^N}U_\varepsilon^2
\log U_\varepsilon^2\,dx
&=(N-2)\varepsilon^2|\log\varepsilon|
\|U_1\|_2^2+\varepsilon^2
\int_{\mathbb R^N}U_1^2\log U_1^2\,dx.
\end{aligned}
\end{align}
For $N \ge 5$, we also have
\begin{align}\label{4.30}
-\frac12\int_{\mathbb R^N}
U_\varepsilon^2\log U_\varepsilon^2\,dx
=
-\frac{N-2}{2}\|U_1\|_2^2
\varepsilon^2|\log\varepsilon|
+O(\varepsilon^2).
\end{align}
\end{lemma}

\begin{proof}
The argument for all norms except the logarithmic term can be found in \cite{ChenJinRadulescuZhou2026}. For the reader's convenience, we give the full proof here. We start with the gradient estimate in $N = 3, 4$. By the fact that the function $U_\varepsilon$ in \eqref{4.26} is an optimizer of the HLS inequality, we have 
\begin{align}\label{4.31}
\int_{\mathbb R^N}\left|\nabla\left[\kappa_{N,\beta}\left(\frac{\varepsilon}{\varepsilon^2+|x|^2}\right)^{\frac{N-2}{2}}\right]\right|^2dx = S_\beta^{\frac{N+\beta}{\beta+2}}.
\end{align}
Specifically, for $N \ge 5$, \eqref{4.31} exactly holds. In the region $|x| < 1$, the cutoff function is exactly \eqref{4.26}. For the region $|x| > 1$, we can write
\begin{align}\label{4.32}
\begin{aligned}
&\int_{|x|\ge1}\left|\nabla\left[\kappa_{N,\beta}\left(\frac{\varepsilon}{\varepsilon^2+|x|^2}\right)^{\frac{N-2}{2}}\right]\right|^2dx\\
&\qquad = \kappa_{N,\beta}^2(N-2)^2\varepsilon^{N-2}\int_{|x|\ge1}\frac{|x|^2}{(\varepsilon^2+|x|^2)^N}dx \\
&\qquad \le \kappa_{N,\beta}^2(N-2)^2|\mathbb S^{N-1}|\varepsilon^{N-2}\int_1^\infty r^{1-N}\,dr\\
&\qquad = O(\varepsilon^{N-2}).
\end{aligned}
\end{align}
In the region $1 < |x| < 2$, we have
\begin{align}\label{4.33}
\begin{aligned}
\int_{1<|x|<2}|\nabla U_\varepsilon|^2dx=\kappa_{N,\beta}^2\left(\frac{\varepsilon}{1+\varepsilon^2}\right)^{N-2}|B_2\setminus B_1| =O(\varepsilon^{N-2}).
\end{aligned}
\end{align}
From \eqref{4.31}, \eqref{4.32}, and \eqref{4.33}, we obtain
\[
\|\nabla U_\varepsilon\|_2^2 = S_\beta^{\frac{N+\beta}{\beta+2}}+O(\varepsilon^{N-2}),
\]
which are the first estimates of the relations \eqref{4.27} and \eqref{4.28} (respectively for $N = 3$ and $N = 4$). For $N \ge 5$, we directly use \eqref{4.31}. Now we estimate the critical Choquard term. For $N \ge 5$, direct calculation on the full bubble gives
\begin{align}\label{4.34}
B_{\beta, q}\left[\kappa_{N,\beta}\left(\frac{\varepsilon}{\varepsilon^2+|x|^2}\right)^{\frac{N-2}{2}}\right] = S_\beta^{\frac{N+\beta}{\beta+2}}.
\end{align}
For $N = 3, 4$, we apply the HLS inequality. Note that $q\frac{2N}{N + \beta} = \frac{2N}{N - 2}$. For $v, w\in D^{1, 2}(\mathbb R^N)$, we have
\begin{align}\label{4.35}
\left|B_{\beta,q}[v]-B_{\beta,q}[w]\right|&\le C_{N,\beta}\left\||v|^q-|w|^q\right\|_{\frac{2N}{N+\beta}}\left(\left\||v|^q\right\|_{\frac{2N}{N+\beta}}+\left\||w|^q\right\|_{\frac{2N}{N+\beta}}\right).
\end{align}
For $|x| \ge 1$, direct calculation gives
\begin{align}\label{4.36}
\begin{aligned}
\int_{|x|\ge1}U_\varepsilon^{\frac{2N}{N-2}}\,dx&=\int_{|x|\ge1}\left[\kappa_{N,\beta}\left(\frac{\varepsilon}{\varepsilon^2+|x|^2}\right)^{\frac{N-2}{2}}\right]^{\frac{2N}{N-2}}dx\\
&=\kappa_{N,\beta}^{\frac{2N}{N-2}}|\mathbb S^{N-1}|\varepsilon^N\int_1^\infty\frac{r^{N-1}}{(\varepsilon^2+r^2)^N}\,dr\\
&\le C\varepsilon^N\int_1^\infty r^{-N-1}\,dr\\
&=O(\varepsilon^N).
\end{aligned}
\end{align}
For $1 < |x| < 2$, we also have
\begin{align}\label{4.37}
\begin{aligned}
\int_{1 < |x| < 2}U_\varepsilon^{\frac{2N}{N-2}}\,dx&=\int_{1 < |x| < 2}\left[\kappa_{N,\beta}\left(\frac{\varepsilon}{\varepsilon^2+|x|^2}\right)^{\frac{N-2}{2}}\right]^{\frac{2N}{N-2}}dx\\
&\le C\varepsilon^N\int_1^2(2-r)^{\frac{2N}{N-2}}r^{N-1}\,dr\\
&=O(\varepsilon^N).
\end{aligned}
\end{align}
Applying \eqref{4.36}-\eqref{4.37} and the fact that for $N = 3, 4$ the bubbles are the same (full bubble), we obtain
\[
\left\|
\left[
\kappa_{N,\beta}
\left(
\frac{\varepsilon}{\varepsilon^2+|x|^2}
\right)^{\frac{N-2}{2}}
\right]^q-U_\varepsilon^q
\right\|_{\frac{2N}{N+\beta}}
=
O\left(\varepsilon^{\frac{N+\beta}{2}}\right).
\]
Applying \eqref{4.35}, we obtain
\[
\begin{aligned}
\left|B_{\beta,q}[U_\varepsilon]-S_\beta^{\frac{N+\beta}{\beta+2}}\right|
\le C\left\|\left[\kappa_{N,\beta}\left(\frac{\varepsilon}{\varepsilon^2+|x|^2}\right)^{\frac{N-2}{2}}\right]^q-U_\varepsilon^q\right\|_{\frac{2N}{N+\beta}}\le C\varepsilon^{\frac{N+\beta}{2}}.
\end{aligned}
\]
Therefore
\[
B_{\beta,q}[U_\varepsilon]= S_\beta^{\frac{N+\beta}{\beta+2}} + O\left(\varepsilon^{\frac{N+\beta}{2}}\right).
\]
It remains to estimate the $L^2$ norm and logarithmic energy. For $N = 3$, using the polar coordinates and the definition of $U_\varepsilon$, we obtain
\begin{align*}
\int_{|x|<1}U_\varepsilon^2\,dx
&=\kappa_{N,\beta}^2|\mathbb S^2|\varepsilon
  \int_0^1\frac{r^2}{\varepsilon^2+r^2}\,dr\\
&=\kappa_{N,\beta}^2|\mathbb S^2|\varepsilon
  \left(1-\varepsilon\arctan\frac1\varepsilon\right)\\
&=\kappa_{N,\beta}^2|\mathbb S^2|\varepsilon
  +O(\varepsilon^2).
\end{align*}
On the annulus $1<|x|<2$, we have
\begin{align*}
\int_{1<|x|<2}U_\varepsilon^2\,dx
&=\frac{\kappa_{N,\beta}^2\varepsilon}{1+\varepsilon^2}
  |\mathbb S^2|\int_1^2r^2(2-r)^2\,dr\\
&=\frac8{15}\kappa_{N,\beta}^2|\mathbb S^2|\varepsilon
  +O(\varepsilon^3).
\end{align*}
Consequently,
\[
\|U_\varepsilon\|_2^2
=\frac{23}{15}\kappa_{N,\beta}^2|\mathbb S^2|
 \varepsilon+O(\varepsilon^2).
\]
For $|x|<1$, we have
\[
\log U_\varepsilon^2(x)
=\log\varepsilon+
 \log\frac{\kappa_{N,\beta}^2}{\varepsilon^2+|x|^2}.
\]
Thus,
\begin{align*}
&\int_{|x|<1}U_\varepsilon^2\log U_\varepsilon^2\,dx
=\kappa_{N,\beta}^2|\mathbb S^2|\varepsilon
 \log\varepsilon\int_0^1\frac{r^2}{\varepsilon^2+r^2}\,dr \\
&\qquad+\kappa_{N,\beta}^2|\mathbb S^2|\varepsilon
 \int_0^1\frac{r^2}{\varepsilon^2+r^2}
 \log\frac{\kappa_{N,\beta}^2}{\varepsilon^2+r^2}\,dr.
\end{align*}
The last integral is uniformly bounded. Indeed, splitting its absolute value into $(0,\varepsilon)$ and $(\varepsilon,1)$ gives
\begin{align*}
&\int_0^\varepsilon\frac{r^2}{\varepsilon^2+r^2}
  \left|\log\frac{\kappa_{N,\beta}^2}{\varepsilon^2+r^2}\right|dr
  \leq C\varepsilon(1+|\log\varepsilon|)\leq C,\\
&\int_\varepsilon^1\frac{r^2}{\varepsilon^2+r^2}
  \left|\log\frac{\kappa_{N,\beta}^2}{\varepsilon^2+r^2}\right|dr
  \leq C\int_0^1(1+|\log r|)\,dr\leq C.
\end{align*}
Therefore,
\[
\int_{|x|<1}U_\varepsilon^2\log U_\varepsilon^2\,dx
=\kappa_{N,\beta}^2|\mathbb S^2|
 \varepsilon\log\varepsilon+O(\varepsilon).
\]
On the annulus, we use
\[
\log U_\varepsilon^2(x)
=\log\varepsilon+
 \log\frac{\kappa_{N,\beta}^2(2-|x|)^2}{1+\varepsilon^2},
\]
and we obtain
\begin{align*}
&\int_{1<|x|<2}U_\varepsilon^2\log U_\varepsilon^2\,dx\\
&\quad=\frac{\kappa_{N,\beta}^2\varepsilon}{1+\varepsilon^2}
 |\mathbb S^2|\int_1^2r^2(2-r)^2
 \left[\log\varepsilon+
 \log\frac{\kappa_{N,\beta}^2(2-r)^2}{1+\varepsilon^2}\right]dr\\
&\quad=\frac8{15}\kappa_{N,\beta}^2|\mathbb S^2|
 \varepsilon\log\varepsilon+O(\varepsilon),
\end{align*}
since $\int_1^2r^2(2-r)^2|\log(2-r)|\,dr<\infty$. Combining the core and annulus contributions gives
\[
\int_{\mathbb R^3}U_\varepsilon^2\log U_\varepsilon^2\,dx
=\frac{23}{15}\kappa_{N,\beta}^2|\mathbb S^2|
 \varepsilon\log\varepsilon+O(\varepsilon).
\]
For $N = 4$, we first compute the $L^2$ contribution from $|x|<1$:
\begin{align*}
\int_{|x|<1}U_\varepsilon^2\,dx
&=\kappa_{N,\beta}^2|\mathbb S^3|\varepsilon^2
  \int_0^1\frac{r^3}{(\varepsilon^2+r^2)^2}\,dr\\
&=\frac{\kappa_{N,\beta}^2|\mathbb S^3|\varepsilon^2}{2}
 \left[
 \log\frac{1+\varepsilon^2}{\varepsilon^2}
 +\frac{\varepsilon^2}{1+\varepsilon^2}-1
 \right]\\
&=\kappa_{N,\beta}^2|\mathbb S^3|\varepsilon^2
  |\log\varepsilon|+O(\varepsilon^2).
\end{align*}
Since
\[
\int_{1<|x|<2}U_\varepsilon^2\,dx
=\frac{\kappa_{N,\beta}^2\varepsilon^2}{(1+\varepsilon^2)^2}
 |\mathbb S^3|\int_1^2r^3(2-r)^2\,dr
=O(\varepsilon^2),
\]
it follows that
\[
\|U_\varepsilon\|_2^2
=\kappa_{N,\beta}^2|\mathbb S^3|
\varepsilon^2|\log\varepsilon|+O(\varepsilon^2).
\]
For $|x|<1$, we have
\[
\log U_\varepsilon^2(x)
=\log\kappa_{N,\beta}^2+2\log\varepsilon
 -2\log(\varepsilon^2+|x|^2).
\]
Consequently,
\begin{align*}
& \int_{|x|<1}U_\varepsilon^2\log U_\varepsilon^2\,dx
=\kappa_{N,\beta}^2|\mathbb S^3|\varepsilon^2 \\
& \qquad \times \Bigg[
 (\log\kappa_{N,\beta}^2+2\log\varepsilon)
 \int_0^1\frac{r^3}{(\varepsilon^2+r^2)^2}\,dr
 -2\int_0^1
 \frac{r^3\log(\varepsilon^2+r^2)}{(\varepsilon^2+r^2)^2}\,dr
 \Bigg].
\end{align*}
The first radial integral satisfies
\[
\int_0^1\frac{r^3}{(\varepsilon^2+r^2)^2}\,dr
=|\log\varepsilon|-\frac12+O(\varepsilon^2).
\]
For the second one, setting $z=\varepsilon^2+r^2$ yields
\begin{align*}
\int_0^1\frac{r^3\log(\varepsilon^2+r^2)}{(\varepsilon^2+r^2)^2}\,dr& =\frac12\int_{\varepsilon^2}^{1+\varepsilon^2}\frac{(t-\varepsilon^2)\log t}{t^2}\,dt\\
&=\frac14\left[\log^2(1+\varepsilon^2)-\log^2(\varepsilon^2)\right]+\frac{\varepsilon^2}{2}\left[\frac{\log(1+\varepsilon^2)+1}{1+\varepsilon^2}-\frac{\log(\varepsilon^2)+1}{\varepsilon^2}\right]\\
&=-|\log\varepsilon|^2+|\log\varepsilon|-\frac12+O(\varepsilon^2).
\end{align*}
Substitution gives
\begin{align*}
\int_{|x|<1}U_\varepsilon^2\log U_\varepsilon^2\,dx&=\kappa_{N,\beta}^2|\mathbb S^3|\varepsilon^2\times \Bigg[(\log\kappa_{N,\beta}^2-2|\log\varepsilon|)\bigg(|\log\varepsilon|-\frac12+O(\varepsilon^2)\bigg)\\
&\qquad -2\bigg(-|\log\varepsilon|^2+|\log\varepsilon|-\frac12+O(\varepsilon^2)\bigg)\Bigg]\\
&=\kappa_{N,\beta}^2|\mathbb S^3|
 (\log\kappa_{N,\beta}^2-1)
 \varepsilon^2|\log\varepsilon|+O(\varepsilon^2).
\end{align*}
In particular, the terms of order $\varepsilon^2|\log\varepsilon|^2$ cancel exactly. For the annular contribution, we have
\begin{align*}
\int_{1<|x|<2}U_\varepsilon^2\log U_\varepsilon^2\,dx
=\frac{\kappa_{N,\beta}^2\varepsilon^2}{(1+\varepsilon^2)^2}
 |\mathbb S^3|\int_1^2r^3(2-r)^2\times\left[
 2\log\varepsilon+
 \log\frac{\kappa_{N,\beta}^2(2-r)^2}{(1+\varepsilon^2)^2}
 \right]dr.
\end{align*}
We get
\[
\int_{1<|x|<2}U_\varepsilon^2\log U_\varepsilon^2\,dx
=-\frac{7}{5}\kappa_{N,\beta}^2|\mathbb S^3|
 \varepsilon^2|\log\varepsilon|+O(\varepsilon^2).
\]
Combining this with the core contribution, we obtain
\[
\int_{\mathbb R^4}U_\varepsilon^2\log U_\varepsilon^2\,dx
=\kappa_{N,\beta}^2|\mathbb S^3|
 \left(\log\kappa_{N,\beta}^2-\frac{12}{5}\right)
 \varepsilon^2|\log\varepsilon|+O(\varepsilon^2).
\]
For $N \ge 5$, the full bubble belongs to $L^2(\mathbb R^N)$. Since $U_\varepsilon(\varepsilon y)=\varepsilon^{-(N-2)/2}U_1(y)$, a change of variables $x=\varepsilon y$ gives
\begin{align*}
\|U_\varepsilon\|_2^2
=\int_{\mathbb R^N}
  \varepsilon^{-(N-2)}U_1^2(x/\varepsilon)\,dx
=\varepsilon^2\int_{\mathbb R^N}U_1^2(y)\,dy
=\varepsilon^2\|U_1\|_2^2.
\end{align*}
Moreover,
\begin{align*}
\int_{\mathbb R^N}U_\varepsilon^2\log U_\varepsilon^2\,dx&=\varepsilon^2\int_{\mathbb R^N}U_1^2(y)
 \log\left(\varepsilon^{-(N-2)}U_1^2(y)\right)\,dy\\
&=-(N-2)\varepsilon^2\log\varepsilon\,
 \|U_1\|_2^2
 +\varepsilon^2\int_{\mathbb R^N}U_1^2\log U_1^2\,dx\\
&=(N-2)\varepsilon^2|\log\varepsilon|\,
 \|U_1\|_2^2
 +\varepsilon^2\int_{\mathbb R^N}U_1^2\log U_1^2\,dx.
\end{align*}
The final integral is finite because $U_1^2(x)\sim\kappa_{N,\beta}^2|x|^{-2(N-2)}$ as $|x|\to\infty$ and $N\geq5$. In particular,
\[
-\frac12\int_{\mathbb R^N}U_\varepsilon^2\log U_\varepsilon^2\,dx = -\frac{N - 2}{2}\varepsilon^2|\log\varepsilon|\, \|U_1\|_2^2+O(\varepsilon^2).
\]
This completes the proof.
\end{proof}

\begin{lemma}\label{lem 4.10}
Assume that $K > \frac{q}{q - 1}$ and $0 < c < \min\{c_F, c_G, c_K\}$. Then
\[
m^-(c) < m^+(c) + \frac{\beta + 2}{2(N + \beta)}S_\beta^{\frac{N + \beta}{\beta + 2}}.
\]
\end{lemma}

\begin{proof}[Proof of Lemma~\ref{lem 4.10} for $N = 3, 4$:]
Let $u_c^+ \in \mathcal{P}_c^+$ be the solution of problem \eqref{eq main}. We have
\[
-\Delta u_c^++\lambda u_c^+ = u_c^+\log (u_c^+)^2 +\mu(I_\alpha*(u_c^+)^p)(u_c^+)^{p - 1} +(I_\beta*(u_c^+)^q)(u_c^+)^{q - 1}.
\]
Testing this equation with $U_\varepsilon$, we obtain
\begin{align}\label{4.38}
\begin{aligned}
&\int_{\mathbb R^N}\nabla u_c^+\cdot\nabla U_\varepsilon\,dx
-\int_{\mathbb R^N}
u_c^+\bigl(\log(u_c^+)^2+1\bigr)U_\varepsilon\,dx\\
&\quad
-\mu\int_{\mathbb R^N}
(I_\alpha*(u_c^+)^p)(u_c^+)^{p-1}U_\varepsilon\,dx\\
&\quad
-\int_{\mathbb R^N}
(I_\beta*(u_c^+)^q)(u_c^+)^{q-1}U_\varepsilon\,dx\\
&=-(\lambda_c^++1)
\int_{\mathbb R^N}u_c^+U_\varepsilon\,dx.
\end{aligned}
\end{align}
By the definition of $U_\varepsilon$ in \eqref{4.25}, we obtain
\begin{align}\label{4.39}
\begin{aligned}
\int_{\mathbb R^N}u_c^+U_\varepsilon\,dx
&\le
\|u_c^+\|_{L^\infty(B_2)}
\int_{B_2}U_\varepsilon\,dx\\
&\le
C\varepsilon^{\frac{N-2}{2}}
\int_0^1
\frac{r^{N-1}}{(\varepsilon^2+r^2)^{\frac{N-2}{2}}}\,dr
+C\varepsilon^{\frac{N-2}{2}}
\int_1^2(2-r)r^{N-1}\,dr\\
&\le
C\varepsilon^{\frac{N-2}{2}}
\int_0^1
r^{N-1-(N-2)}\,dr
+C\varepsilon^{\frac{N-2}{2}}\\
&\le C\varepsilon^{\frac{N-2}{2}},
\end{aligned}
\end{align}
where we have used $(\varepsilon^2+r^2)^{\frac{N-2}{2}}\ge r^{N-2}$. By Lemma~\ref{lem 4.9}, we also have $\|U_\varepsilon\|_2^2=o\!\left(\varepsilon^{\frac{N-2}{2}}\right)$. Therefore for $t \in [0, T]$ with $T < +\infty$, we obtain
\begin{align}\label{4.40}
\begin{aligned}
\|u_c^++tU_\varepsilon\|_2^2
&=\|u_c^+\|_2^2
+2t\int_{\mathbb R^N}u_c^+U_\varepsilon\,dx
+t^2\|U_\varepsilon\|_2^2\\
&=c^2+
2t\int_{\mathbb R^N}u_c^+U_\varepsilon\,dx
+o\!\left(\varepsilon^{\frac{N-2}{2}}\right).
\end{aligned}
\end{align}
For $a, b\ge 0$ and $q \ge 2$, we use the basic inequalities
\[
(a+b)^q-a^q-b^q\ge qa^{q-1}b,
\qquad
(a+b)^q-a^q-b^q\ge qb^{q-1}a,
\]
and we obtain
\begin{align}\label{4.41}
\begin{aligned}
B_{\beta,q}[u_c^++tU_\varepsilon]
&=
\int_{\mathbb R^N}
\bigl(I_\beta*(u_c^++tU_\varepsilon)^q\bigr)
(u_c^++tU_\varepsilon)^q\,dx\\
&\ge
B_{\beta,q}[u_c^+]+t^{2q}B_{\beta,q}[U_\varepsilon]\\
&\quad+
2\int_{\mathbb R^N}
(I_\beta*(u_c^+)^q)
\bigl[(u_c^++tU_\varepsilon)^q
-(u_c^+)^q-t^qU_\varepsilon^q\bigr]\,dx\\
&\quad+
2t^q\int_{\mathbb R^N}
(I_\beta*U_\varepsilon^q)
\bigl[(u_c^++tU_\varepsilon)^q
-(u_c^+)^q-t^qU_\varepsilon^q\bigr]\,dx\\
&\ge
B_{\beta,q}[u_c^+]+t^{2q}B_{\beta,q}[U_\varepsilon]\\
&\quad+
2qt\int_{\mathbb R^N}
(I_\beta*(u_c^+)^q)(u_c^+)^{q-1}U_\varepsilon\,dx\\
&\quad+
2qt^{2q-1}\int_{\mathbb R^N}
(I_\beta*U_\varepsilon^q)
U_\varepsilon^{q-1}u_c^+\,dx.
\end{aligned}
\end{align}
Using the positivity and continuity of $u_c^+$, we get
\begin{align}\label{4.42}
\begin{aligned}
&\int_{\mathbb R^N}
(I_\beta*U_\varepsilon^q)
U_\varepsilon^{q-1}u_c^+\,dx\\
&\ge
\left(\min_{|x|\le1}u_c^+(x)\right)
\int_{|x|\le\varepsilon}
\int_{|y|\le\varepsilon}
I_\beta(x-y)
U_\varepsilon(y)^q
U_\varepsilon(x)^{q-1}\,dy\,dx\\
&=
C\varepsilon^{2N}
\varepsilon^{\beta-N}
\varepsilon^{-\frac{N+\beta}{2}}
\varepsilon^{-\frac{\beta+2}{2}}
\int_{B_1}\int_{B_1}
\frac{
dx\,dy
}{
|x-y|^{N-\beta}
(1+|y|^2)^{\frac{N+\beta}{2}}
(1+|x|^2)^{\frac{\beta+2}{2}}
}\\
&\ge C\varepsilon^{\frac{N-2}{2}}.
\end{aligned}
\end{align}
Applying \eqref{4.42} on \eqref{4.41}, we obtain
\begin{align}\label{4.43}
-t^{2q-1}
\int_{\mathbb R^N}
(I_\beta*U_\varepsilon^q)
U_\varepsilon^{q-1}u_c^+\,dx
\le
-Ct^{2q-1}\varepsilon^{\frac{N-2}{2}}.
\end{align}
Continuing, applying \eqref{4.40} and using $\|U_\varepsilon\|_2^2=o\left(\varepsilon^{\frac{N-2}{2}}\right)$ we obtain
\begin{align}\label{4.44}
\begin{aligned}
\frac{c^2}{\|u_c^++tU_\varepsilon\|_2^2}
&=
1-\frac{2t}{c^2}
\int_{\mathbb R^N}u_c^+U_\varepsilon\,dx
+o\left(\varepsilon^{\frac{N-2}{2}}\right),\\[2mm]
\log\frac{\|u_c^++tU_\varepsilon\|_2}{c}
&=
\frac{t}{c^2}
\int_{\mathbb R^N}u_c^+U_\varepsilon\,dx
+o\left(\varepsilon^{\frac{N-2}{2}}\right),\\[2mm]
\left(
\frac{\|u_c^++tU_\varepsilon\|_2}{c}
\right)^{(N-2)p-(N+\alpha)}
&=
1-
\frac{N+\alpha-(N-2)p}{c^2}\,
t\int_{\mathbb R^N}u_c^+U_\varepsilon\,dx
+o\left(\varepsilon^{\frac{N-2}{2}}\right),
\end{aligned}
\end{align}
uniformly for $t \in [0, 2]$. Now, define
\[
z_{c, t}:= \theta_{c, t}^{\frac{N - 2}{2}}(u_c^++tU_\varepsilon)\left(\theta_{c, t}\,x\right), \qquad \theta_{c, t}:= \frac{\|u_c^++tU_\varepsilon\|_2}{c},
\]
and notice that $z_{c, t} \in S_c$. From \eqref{energy}, we can write
\begin{align}\label{4.45}
\begin{aligned}
E(z_{c, t})&=
\frac12\|\nabla(u_c^++tU_\varepsilon)\|_2^2
+\frac{c^2}{2}\\
&-\frac{c^2}{2\|u_c^++tU_\varepsilon\|_2^2}
\int_{\mathbb R^N}
(u_c^++tU_\varepsilon)^2
\log (u_c^++tU_\varepsilon)^2\,dx\\
&-\frac{N-2}{2}c^2
\log\frac{\|u_c^++tU_\varepsilon\|_2}{c}\\
&-\frac{\mu}{2p}
\left(
\frac{\|u_c^++tU_\varepsilon\|_2}{c}
\right)^{(N-2)p-(N+\alpha)}
B_{\alpha,p}[u_c^++tU_\varepsilon]\\
&-\frac1{2q}B_{\beta,q}[u_c^++tU_\varepsilon].
\end{aligned}
\end{align}
For the logarithmic term in \eqref{4.45}, we apply the elementary identity
\[
(a+b)^2\log(a+b)^2\ge a^2\log a^2+2ab(\log a^2+1)+b^2\log b^2,
\]
and we obtain
\begin{align}\label{4.46}
\begin{aligned}
&\int_{\mathbb R^N}
(u_c^++tU_\varepsilon)^2
\log(u_c^++tU_\varepsilon)^2\,dx\\
&\ge
\int_{\mathbb R^N}(u_c^+)^2\log(u_c^+)^2\,dx
+2t\int_{\mathbb R^N}
u_c^+\bigl(\log(u_c^+)^2+1\bigr)U_\varepsilon\,dx\\
&\qquad+
t^2\int_{\mathbb R^N}
U_\varepsilon^2\log U_\varepsilon^2\,dx
+t^2\log(t^2)\|U_\varepsilon\|_2^2,
\end{aligned}
\end{align}
and for the subcritical Choquard term, we have
\begin{align}\label{4.47}
\begin{aligned}
B_{\alpha,p}[u_c^++tU_\varepsilon]
&\ge
B_{\alpha,p}[u_c^+]
+2pt\int_{\mathbb R^N}
(I_\alpha*(u_c^+)^p)(u_c^+)^{p-1}U_\varepsilon\,dx.
\end{aligned}
\end{align}
Applying \eqref{4.43}, \eqref{4.44}, \eqref{4.46}, and \eqref{4.47} on \eqref{4.45}, we obtain
\begin{align}\label{4.48}
\begin{aligned}
E(z_{c, t})&\le
m^+(c)+\frac{t^2}{2}\|\nabla U_\varepsilon\|_2^2-\frac{t^{2q}}{2q}B_{\beta,q}[U_\varepsilon]\\
&+t\int_{\mathbb R^N}\nabla u_c^+\cdot\nabla U_\varepsilon\,dx-t\int_{\mathbb R^N}
u_c^+\bigl(\log(u_c^+)^2+1\bigr)U_\varepsilon\,dx\\
&-\mu t\int_{\mathbb R^N}
(I_\alpha*(u_c^+)^p)(u_c^+)^{p-1}U_\varepsilon\,dx\\
&-t\int_{\mathbb R^N}
(I_\beta*(u_c^+)^q)(u_c^+)^{q-1}U_\varepsilon\,dx\\
&-t^{2q-1}
\int_{\mathbb R^N}
(I_\beta*U_\varepsilon^q)U_\varepsilon^{q-1}u_c^+\,dx\\
&+\frac{t}{c^2}
\left(\int_{\mathbb R^N}u_c^+U_\varepsilon\,dx\right)
\Bigg[
\int_{\mathbb R^N}(u_c^+)^2\log(u_c^+)^2\,dx-\frac{N-2}{2}c^2\\
&\quad+\frac{\mu[N+\alpha-(N-2)p]}{2p}B_{\alpha,p}[u_c^+]
\Bigg]\\
&-\frac{t^2}{2}\int_{\mathbb R^N}
U_\varepsilon^2\log U_\varepsilon^2\,dx
-\frac{t^2\log(t^2)}{2}\|U_\varepsilon\|_2^2
+o\left(\varepsilon^{\frac{N-2}{2}}\right),
\end{aligned}
\end{align}
uniformly for $t \in [0, 2]$. Now, we simplify \eqref{4.48}. From the main equation, we can directly obtain the system
\[
\begin{cases}
\displaystyle{\|\nabla u_c^+\|_2^2+\lambda_c^+c^2 = \int_{\mathbb R^N}(u_c^+)^2\log(u_c^+)^2\,dx +\mu B_{\alpha,p}[u_c^+] +B_{\beta,q}[u_c^+]}, \\
\displaystyle{\|\nabla u_c^+\|_2^2 = \frac N2c^2+\mu\eta_pB_{\alpha,p}[u_c^+] +B_{\beta,q}[u_c^+]},
\end{cases}
\]
and from which we get
\begin{align}\label{4.49}
\begin{aligned}
&\int_{\mathbb R^N}(u_c^+)^2\log(u_c^+)^2\,dx
-\frac{N-2}{2}c^2
+\frac{\mu[N+\alpha-(N-2)p]}{2p}B_{\alpha,p}[u_c^+]\\
&\qquad=(\lambda_c^++1)c^2.
\end{aligned}
\end{align}
Combining \eqref{4.38} and \eqref{4.49} and simplifying, all the terms with $\displaystyle\int_{\mathbb R^N}u_c^+U_\varepsilon\,dx$ cancel and we get
\begin{align}\label{4.50}
\begin{aligned}
E(z_{c, t})
&\le
m^+(c)
+\frac{t^2}{2}\|\nabla U_\varepsilon\|_2^2
-\frac{t^{2q}}{2q}B_{\beta,q}[U_\varepsilon]\\
&-t^{2q-1}\int_{\mathbb R^N}
(I_\beta*U_\varepsilon^q)U_\varepsilon^{q-1}u_c^+\,dx\\
&-\frac{t^2}{2}
\int_{\mathbb R^N}U_\varepsilon^2\log U_\varepsilon^2\,dx
-\frac{t^2\log(t^2)}2\|U_\varepsilon\|_2^2\\
&+o\left(\varepsilon^{\frac{N-2}{2}}\right).
\end{aligned}
\end{align}
Applying the estimates \eqref{4.27} and \eqref{4.28} on \eqref{4.50}, we obtain
\[
\begin{aligned}
E(z_{c, t})&\le m^+(c)+S_\beta^{\frac q{q-1}}\left(\frac{t^2}{2}-\frac{t^{2q}}{2q}\right)-Ct^{2q-1}\varepsilon^{\frac{N-2}{2}}+o\left(\varepsilon^{\frac{N-2}{2}}\right).
\end{aligned}
\]
Since $q > 1$, the function $t\mapsto \frac{t^2}{2}-\frac{t^{2q}}{2q}$ attains its unique maximum $\frac{q-1}{2q}$ at $t=1$. For $t\in[\frac12,\frac32]$, we have
\[
\begin{aligned}
E(z_{c,t})
&\le m^+(c)+\frac{q-1}{2q}S_\beta^{\frac q{q-1}}
-C2^{1-2q}\varepsilon^{\frac{N-2}{2}}
+o\left(\varepsilon^{\frac{N-2}{2}}\right)\\
&<m^+(c)+\frac{q-1}{2q}S_\beta^{\frac q{q-1}},
\end{aligned}
\]
for sufficiently small $\varepsilon>0$. On the other hand,
\[
\max_{t\in[0 ,\frac{1}{2}]\cup[\frac{3}{2},2]} \left(\frac{t^2}{2}-\frac{t^{2q}}{2q}\right) <\frac{q-1}{2q}.
\]
Hence, using the uniformity of the remainder, we obtain
\[
\max_{t\in[0,2]}E(z_{c,t}) <m^+(c)+\frac{q-1}{2q}S_\beta^{\frac q{q-1}}.
\]
Moreover, $z_{c,0}=u_c^+$ and $\frac{2^2}{2}-\frac{2^{2q}}{2q}<0$, so $E(z_{c,2})<m^+(c)$ for small $\varepsilon>0$. The path $t\mapsto z_{c,t}$ is continuous in $H_r$.
Since $z_{c,0}=u_c^+\in W_{c,K}$ and $E(z_{c,2})<m^+(c)$, {Lemma~\ref{lem 4.7} implies $z_{c,2} \not \in W_{c,K}$ and so $\|\nabla z_{c,2}\|_2>t_0$.} If $P(z_{c,2})\ge0$, Lemma~\ref{lem 4.5} gives $s^+(z_{c,2})<0\le s^-(z_{c,2})$ and
\[
 E(z_{c,2})\ge E(s^+(z_{c,2})\star z_{c,2})\ge m^+(c),
\]
a contradiction. {Thus $P(z_{c,2})<0$ and $\|\nabla z_{c,2}\|_2 \not =t_0$ (by Lemma~\ref{lem 4.1}).} At the last $t_*\in(0,2)$ with $\|\nabla z_{c,t_*}\|_2=t_0$, Lemma~\ref{lem 4.1} ensures $P(z_{c,t_*})>0$. By continuity, some $t_-\in(t_*,2)$ satisfies $P(z_{c,t_-})=0$ and $\|\nabla z_{c,t_-}\|_2>t_0$; hence $z_{c,t_-}\in\mathcal P_c^-$. Therefore
\[
m^-(c)\le\max_{t\in[0,2]}E(z_{c,t}) <m^+(c)+ \frac{\beta+2}{2(N+\beta)} S_\beta^{\frac{N+\beta}{\beta+2}}.
\]
This completes the proof.
\end{proof}

\begin{proof}[Proof of Lemma~\ref{lem 4.10} for $N \ge 5$:]
For $0\le t\le 2$, we consider
\begin{equation}\label{eq:gap-high-path}
\begin{aligned}
z_{c,t}(x)
&=\left(\frac{\|u_c^++tU_\varepsilon\|_2}{c}\right)^{\frac{N-2}{2}}
(u_c^++tU_\varepsilon)
\left(\frac{\|u_c^++tU_\varepsilon\|_2}{c}\,x\right).
\end{aligned}
\end{equation}
By the choice of the spatial dilation,
\begin{equation}\label{eq:gap-high-invariances}
\begin{aligned}
\|z_{c,t}\|_2&=c,\\
\|\nabla z_{c,t}\|_2^2
&=\|\nabla(u_c^++tU_\varepsilon)\|_2^2,\\
B_{\beta,q}[z_{c,t}]&=B_{\beta,q}[u_c^++tU_\varepsilon].
\end{aligned}
\end{equation}
We first record a lower bound for the Choquard interactions that remains valid for every $r>1$. Indeed, for $a,b\ge0$,
\[
(a+b)^r-a^r\ge b^r,
\qquad
(a+b)^r-a^r\ge ra^{r-1}b.
\]
Using the positivity of the Riesz kernel, we obtain
\begin{equation}\label{eq:gap-high-choquard}
\begin{aligned}
&B_{\sigma,r}[u_c^++tU_\varepsilon]-B_{\sigma,r}[u_c^+]\\
={}&2\int_{\mathbb R^N}(I_\sigma*(u_c^+)^r)
\bigl[(u_c^++tU_\varepsilon)^r-(u_c^+)^r\bigr] \,dx\\
&+\int_{\mathbb R^N}
\Bigl(I_\sigma*\bigl[(u_c^++tU_\varepsilon)^r-(u_c^+)^r\bigr]\Bigr)
\bigl[(u_c^++tU_\varepsilon)^r-(u_c^+)^r\bigr] \,dx\\
\ge{}&2rt\int_{\mathbb R^N}(I_\sigma*(u_c^+)^r)
(u_c^+)^{r-1}U_\varepsilon\,dx
+t^{2r}B_{\sigma,r}[U_\varepsilon].
\end{aligned}
\end{equation}
In particular, \eqref{eq:gap-high-choquard} applies to
$(\sigma,r)=(\alpha,p)$ and $(\sigma,r)=(\beta,q)$,
without imposing $q\ge2$.\\

\noindent
By Lemma~\ref{lem 4.9}, the full HLS bubble satisfies
\begin{equation}\label{eq:gap-high-bubble}
\begin{aligned}
\|\nabla U_\varepsilon\|_2^2
&=B_{\beta,q}[U_\varepsilon]
=S_\beta^{\frac q{q-1}},\\
\|U_\varepsilon\|_2^2
&=\varepsilon^2\|U_1\|_2^2,\\
\int_{\mathbb R^N}U_\varepsilon^2\log U_\varepsilon^2\,dx
&=(N-2)\varepsilon^2|\log\varepsilon|\|U_1\|_2^2
+\varepsilon^2\int_{\mathbb R^N}U_1^2\log U_1^2\,dx.
\end{aligned}
\end{equation}
We next establish the weighted estimate
\begin{equation}\label{eq:gap-high-weighted}
\begin{aligned}
&\int_{\mathbb R^N}u_c^+U_\varepsilon\,dx
+\left|\int_{\mathbb R^N}
 u_c^+\bigl(\log (u_c^+)^2+1\bigr)U_\varepsilon\,dx\right|\\
&\quad+\int_{\mathbb R^N}
(I_\alpha*(u_c^+)^p)(u_c^+)^{p-1}U_\varepsilon\,dx
=O\left(\varepsilon^{\frac{N-2}{2}}\right).
\end{aligned}
\end{equation}
{Using Lemma~\ref{lem 2.6} we can deduce that $u_c^+$ and $I_\alpha*(u_c^+)^p$ are locally bounded. Since \(u_c^+\in H_r^1(\mathbb R^N)\), the Strauss' radial lemma yields \(u_c^+\in L^\infty(\mathbb R^N)\).	
	 Furthermore,
\[
U_\varepsilon(x)
\le C\varepsilon^{\frac{N-2}{2}}|x|^{-(N-2)}
\qquad (x\ne0),
\]
and} the contributions over $B_1$ satisfy
\[
\int_{B_1}U_\varepsilon\,dx
\le C\varepsilon^{\frac{N-2}{2}}
\int_0^1r^{N-1-(N-2)}\,dr
\le C\varepsilon^{\frac{N-2}{2}}.
\]
On $\mathbb R^N\setminus B_1$, the assumption $N>4$ yields
\[
\int_{|x|>1}|x|^{-2(N-2)}\,dx=|\mathbb S^{N-1}|\int_1^\infty r^{3-N}\,dr<\infty.
\]
Hence, by H\"older's inequality,
\begin{equation}\label{eq:gap-high-L2-weight}
\begin{aligned}
\int_{|x|>1}u_c^+U_\varepsilon\,dx
&\le C\varepsilon^{\frac{N-2}{2}}
\int_{|x|>1}u_c^+(x)|x|^{-(N-2)}\,dx\\
&\le C\varepsilon^{\frac{N-2}{2}}\|u_c^+\|_2
\left(\int_{|x|>1}|x|^{-2(N-2)}\,dx\right)^{1/2}\\
&\le C\varepsilon^{\frac{N-2}{2}}.
\end{aligned}
\end{equation}
Moreover, the boundedness of $u_c^+$ implies
\[
\left|u_c^+\bigl(\log (u_c^+)^2+1\bigr)\right|
\le C\left((u_c^+)^{\frac{N+4}{2N}}+u_c^+\right).
\]
Since $4N(N-2)/(3N-4)>N$ for $N>4$, another application of H\"older's inequality gives
\begin{equation}\label{eq:gap-high-log-weight}
\begin{aligned}
&\int_{|x|>1}(u_c^+)^{\frac{N+4}{2N}}U_\varepsilon\,dx\\
&\le C\varepsilon^{\frac{N-2}{2}}
\int_{|x|>1}(u_c^+)^{\frac{N+4}{2N}}|x|^{-(N-2)}\,dx\\
&\le C\varepsilon^{\frac{N-2}{2}}
\|u_c^+\|_2^{\frac{N+4}{2N}}
\left(\int_{|x|>1}|x|^{-\frac{4N(N-2)}{3N-4}}\,dx
\right)^{\frac{3N-4}{4N}}\\
&\le C\varepsilon^{\frac{N-2}{2}}.
\end{aligned}
\end{equation}
Finally, the HLS inequality and H\"older's inequalities yield
\begin{equation}\label{eq:gap-high-subcritical-weight}
\begin{aligned}
&\int_{|x|>1}(I_\alpha*(u_c^+)^p)(u_c^+)^{p-1}
U_\varepsilon\,dx\\
&\le C\varepsilon^{\frac{N-2}{2}}
\int_{|x|>1}(I_\alpha*(u_c^+)^p)(u_c^+)^{p-1}
|x|^{-(N-2)}\,dx\\
&\le C\varepsilon^{\frac{N-2}{2}}
\|I_\alpha*(u_c^+)^p\|_{\frac{2N}{N-\alpha}}
\|u_c^+\|_{\frac{2Np}{N+\alpha}}^{p-1}
\left(\int_{|x|>1}|x|^{-\frac{2Np(N-2)}{N+\alpha}}\,dx
\right)^{\frac{N+\alpha}{2Np}}\\
&\le C\varepsilon^{\frac{N-2}{2}}
\|u_c^+\|_{\frac{2Np}{N+\alpha}}^{2p-1}
\left(\int_1^\infty
r^{N-1-\frac{2Np(N-2)}{N+\alpha}}\,dr
\right)^{\frac{N+\alpha}{2Np}}\\
&\le C\varepsilon^{\frac{N-2}{2}}.
\end{aligned}
\end{equation}
Here the last integral is finite, since $p>\frac{N+\alpha}{N}$, $2(N-2)>N$ and hence $\frac{2Np(N-2)}{N+\alpha}>N$. Also, $u_c^+\in L^{\frac{2Np}{N+\alpha}}(\mathbb R^N)$ under the hypotheses on $p$. Together, \eqref{eq:gap-high-L2-weight}--\eqref{eq:gap-high-subcritical-weight} prove \eqref{eq:gap-high-weighted}. For $t\in[0,2]$, the mass satisfies
\begin{equation}\label{eq:gap-high-mass}
\|u_c^++tU_\varepsilon\|_2^2
=c^2+2t\int_{\mathbb R^N}u_c^+U_\varepsilon\,dx
+t^2\varepsilon^2\|U_1\|_2^2.
\end{equation}
By \eqref{eq:gap-high-weighted},
\[
\left|2t\int_{\mathbb R^N}u_c^+U_\varepsilon\,dx
+t^2\varepsilon^2\|U_1\|_2^2\right|^2
=O(\varepsilon^3)
\]
uniformly for $t\in[0,2]$, because $N\ge5$. Consequently, Taylor's formula gives, uniformly in $t\in[0,2]$,
\begin{equation}\label{eq:gap-high-mass-expansions}
\begin{aligned}
\frac{c^2}{\|u_c^++tU_\varepsilon\|_2^2}
&=1-\frac{2t}{c^2}\int_{\mathbb R^N}u_c^+U_\varepsilon\,dx
-\frac{t^2\varepsilon^2}{c^2}\|U_1\|_2^2
+O(\varepsilon^3),\\
\log\frac{\|u_c^++tU_\varepsilon\|_2}{c}
&=\frac{t}{c^2}\int_{\mathbb R^N}u_c^+U_\varepsilon\,dx
+\frac{t^2\varepsilon^2}{2c^2}\|U_1\|_2^2
+O(\varepsilon^3),\\
\left(\frac{\|u_c^++tU_\varepsilon\|_2}{c}\right)^
{(N-2)p-(N+\alpha)}
&=1-\frac{N+\alpha-(N-2)p}{c^2}\,
 t\int_{\mathbb R^N}u_c^+U_\varepsilon\,dx\\
&\quad-\frac{N+\alpha-(N-2)p}{2c^2}\,
 t^2\varepsilon^2\|U_1\|_2^2
+O(\varepsilon^3).
\end{aligned}
\end{equation}
The definition of $z_{c,t}$ and the scaling properties imply
\begin{equation}\label{eq:gap-high-exact-energy}
\begin{aligned}
E(z_{c,t})
={}&\frac12\|\nabla(u_c^++tU_\varepsilon)\|_2^2
+\frac{c^2}{2}\\
&-\frac{c^2}{2\|u_c^++tU_\varepsilon\|_2^2}
\int_{\mathbb R^N}(u_c^++tU_\varepsilon)^2
\log(u_c^++tU_\varepsilon)^2\,dx\\
&-\frac{N-2}{2}c^2
\log\frac{\|u_c^++tU_\varepsilon\|_2}{c}\\
&-\frac{\mu}{2p}
\left(\frac{\|u_c^++tU_\varepsilon\|_2}{c}\right)^
{(N-2)p-(N+\alpha)}
B_{\alpha,p}[u_c^++tU_\varepsilon]\\
&-\frac1{2q}B_{\beta,q}[u_c^++tU_\varepsilon].
\end{aligned}
\end{equation}
For all $a,b\ge0$, with the convention $s^2\log s^2=0$ at $s=0$,
\[
(a+b)^2\log(a+b)^2
\ge a^2\log a^2+2ab(\log a^2+1)+b^2\log b^2.
\]
It follows that
\begin{equation}\label{eq:gap-high-log-split}
\begin{aligned}
&\int_{\mathbb R^N}(u_c^++tU_\varepsilon)^2
\log(u_c^++tU_\varepsilon)^2\,dx\\
&\ge \int_{\mathbb R^N}(u_c^+)^2\log(u_c^+)^2\,dx
+2t\int_{\mathbb R^N}u_c^+
\bigl(\log(u_c^+)^2+1\bigr)U_\varepsilon\,dx\\
&\quad+t^2\int_{\mathbb R^N}U_\varepsilon^2\log U_\varepsilon^2\,dx
+t^2\log(t^2)\|U_\varepsilon\|_2^2.
\end{aligned}
\end{equation}
The Choquard lower bounds \eqref{eq:gap-high-choquard} give
\begin{equation}\label{eq:gap-high-interaction-splits}
\begin{aligned}
B_{\alpha,p}[u_c^++tU_\varepsilon]
&\ge B_{\alpha,p}[u_c^+]
+2pt\int_{\mathbb R^N}(I_\alpha*(u_c^+)^p)
(u_c^+)^{p-1}U_\varepsilon\,dx,\\
B_{\beta,q}[u_c^++tU_\varepsilon]
&\ge B_{\beta,q}[u_c^+]
+t^{2q}B_{\beta,q}[U_\varepsilon]\\
&\quad+2qt\int_{\mathbb R^N}(I_\beta*(u_c^+)^q)
(u_c^+)^{q-1}U_\varepsilon\,dx.
\end{aligned}
\end{equation}
Inserting \eqref{eq:gap-high-mass-expansions},
\eqref{eq:gap-high-log-split}, and
\eqref{eq:gap-high-interaction-splits} into
\eqref{eq:gap-high-exact-energy}, and using
\eqref{eq:gap-high-weighted}, we obtain, uniformly for $t\in[0,2]$, hence
\begin{equation}\label{eq:gap-high-long-expansion}
\begin{aligned}
E(z_{c,t})
\le{}&m^+(c)
+\frac{t^2}{2}\|\nabla U_\varepsilon\|_2^2
-\frac{t^{2q}}{2q}B_{\beta,q}[U_\varepsilon]\\
&-\frac{t^2}{2}\int_{\mathbb R^N}
U_\varepsilon^2\log U_\varepsilon^2\,dx
-\frac{t^2\log(t^2)}2\|U_\varepsilon\|_2^2\\
&+t\int_{\mathbb R^N}\nabla u_c^+\cdot\nabla U_\varepsilon\,dx
-t\int_{\mathbb R^N}u_c^+
\bigl(\log(u_c^+)^2+1\bigr)U_\varepsilon\,dx\\
&-\mu t\int_{\mathbb R^N}(I_\alpha*(u_c^+)^p)
(u_c^+)^{p-1}U_\varepsilon\,dx
-t\int_{\mathbb R^N}(I_\beta*(u_c^+)^q)
(u_c^+)^{q-1}U_\varepsilon\,dx\\
&+\left(
\frac{t}{c^2}\int_{\mathbb R^N}u_c^+U_\varepsilon\,dx
+\frac{t^2}{2c^2}\|U_\varepsilon\|_2^2\right)\\
&\qquad\times\bigg[
\int_{\mathbb R^N}(u_c^+)^2\log(u_c^+)^2\,dx
-\frac{N-2}{2}c^2
+\frac{\mu[N+\alpha-(N-2)p]}{2p}
B_{\alpha,p}[u_c^+]\bigg]\\
&\qquad+O(\varepsilon^3).
\end{aligned}
\end{equation}
All remainder estimates above are uniform with respect to $t\in[0,2]$. Since $u_c^+$ solves the normalized Euler--Lagrange equation, testing that equation against $U_\varepsilon$ yields
\begin{equation}\label{eq:gap-high-euler}
\begin{aligned}
&\int_{\mathbb R^N}\nabla u_c^+\cdot\nabla U_\varepsilon\,dx
-\int_{\mathbb R^N}u_c^+
\bigl(\log(u_c^+)^2+1\bigr)U_\varepsilon\,dx\\
&\quad-\mu\int_{\mathbb R^N}(I_\alpha*(u_c^+)^p)
(u_c^+)^{p-1}U_\varepsilon\,dx\\
&\quad-\int_{\mathbb R^N}(I_\beta*(u_c^+)^q)
(u_c^+)^{q-1}U_\varepsilon\,dx
=-(\lambda+1)\int_{\mathbb R^N}u_c^+U_\varepsilon\,dx.
\end{aligned}
\end{equation}
Moreover, the Pohozaev identity $P(u_c^+)=0$ reads
\[
\|\nabla u_c^+\|_2^2
=\frac N2c^2+\mu\eta_pB_{\alpha,p}[u_c^+]
+B_{\beta,q}[u_c^+].
\]
Testing the Euler--Lagrange equation against $u_c^+$ gives
\[
\|\nabla u_c^+\|_2^2+\lambda c^2
=\int_{\mathbb R^N}(u_c^+)^2\log(u_c^+)^2\,dx
+\mu B_{\alpha,p}[u_c^+]+B_{\beta,q}[u_c^+].
\]
Using $1-\eta_p=\frac{N+\alpha-(N-2)p}{2p}$, we infer that
\begin{equation}\label{eq:gap-high-bracket}
\begin{aligned}
&\int_{\mathbb R^N}(u_c^+)^2\log(u_c^+)^2\,dx
-\frac{N-2}{2}c^2\\
&\qquad+\frac{\mu[N+\alpha-(N-2)p]}{2p}
B_{\alpha,p}[u_c^+]
=(\lambda+1)c^2.
\end{aligned}
\end{equation}
Substituting \eqref{eq:gap-high-euler} and \eqref{eq:gap-high-bracket} into \eqref{eq:gap-high-long-expansion}, all terms linear in $\displaystyle\int_{\mathbb R^N}u_c^+U_\varepsilon\,dx$ cancel exactly. Consequently,
\begin{equation}\label{eq:gap-high-postcancel}
\begin{aligned}
E(z_{c,t})
\le{}&m^+(c)
+\frac{t^2}{2}\|\nabla U_\varepsilon\|_2^2
-\frac{t^{2q}}{2q}B_{\beta,q}[U_\varepsilon]\\
&-\frac{t^2}{2}\int_{\mathbb R^N}
U_\varepsilon^2\log U_\varepsilon^2\,dx
+\frac{t^2}{2}\bigl(\lambda+1-\log(t^2)\bigr)
\|U_\varepsilon\|_2^2
+O(\varepsilon^3).
\end{aligned}
\end{equation}
Combining this inequality with \eqref{eq:gap-high-bubble}, we obtain
\begin{equation}\label{eq:gap-high-sharp-energy}
\begin{aligned}
E(z_{c,t})
\le{}&m^+(c)
+S_\beta^{\frac q{q-1}}
\left(\frac{t^2}{2}-\frac{t^{2q}}{2q}\right)\\
&-\frac{N-2}{2}t^2\|U_1\|_2^2
\varepsilon^2|\log\varepsilon|\\
&+\frac{t^2\varepsilon^2}{2}
\left[\bigl(\lambda+1-\log(t^2)\bigr)\|U_1\|_2^2
-\int_{\mathbb R^N}U_1^2\log U_1^2\,dx\right]
+O(\varepsilon^3).
\end{aligned}
\end{equation}
With the convention $t^2\log(t^2)=0$ at $t=0$, the $\varepsilon^2$-coefficient is bounded for $t\in[0,2]$. Hence
\begin{equation}\label{eq:gap-high-uniform-energy}
\begin{aligned}
E(z_{c,t})
\le{}&m^+(c)
+S_\beta^{\frac q{q-1}}
\left(\frac{t^2}{2}-\frac{t^{2q}}{2q}\right)\\
&-\frac{N-2}{2}t^2\|U_1\|_2^2
\varepsilon^2|\log\varepsilon|
+O(\varepsilon^2),
\end{aligned}
\end{equation}
uniformly for $t\in[0, 2]$. For every $q>1$, the function $t\mapsto \frac{t^2}{2}-\frac{t^{2q}}{2q}$ has its unique global maximum $\frac{q-1}{2q}$ at $t=1$. If $t\in[\frac{1}{2}, \frac{3}{2}]$, \eqref{eq:gap-high-uniform-energy} gives
\[
E(z_{c,t})
\le m^+(c)+\frac{q-1}{2q}S_\beta^{\frac q{q-1}}
-\frac{N-2}{8}\|U_1\|_2^2
\varepsilon^2|\log\varepsilon|+O(\varepsilon^2)
<m^+(c)+\frac{q-1}{2q}S_\beta^{\frac q{q-1}}
\]
for all sufficiently small $\varepsilon>0$. On the complementary set $[0, \frac{1}{2}]\cup[\frac{3}{2}, 2]$, one has
\[
\max_{t\in[0,1/2]\cup[3/2,2]}\left(\frac{t^2}{2}-\frac{t^{2q}}{2q}\right)<\frac{q-1}{2q}.
\]
The uniform remainder in \eqref{eq:gap-high-uniform-energy} then implies
\begin{equation}\label{eq:gap-high-max}
\max_{t\in[0,2]}E(z_{c,t})<m^+(c)+\frac{q-1}{2q}S_\beta^{\frac q{q-1}}
\end{equation}
for sufficiently small $\varepsilon>0$. Finally, $z_{c,0}=u_c^+$, and, since $q>1$, then $\frac{2^2}{2}-\frac{2^{2q}}{2q}<0$. It follows from \eqref{eq:gap-high-uniform-energy} that $E(z_{c,2})<m^+(c)$ for all sufficiently small $\varepsilon>0$. The path $t\mapsto z_{c,t}$ is continuous in $H_r$ by continuity of addition and mass normalization in $H_r$.
Since $z_{c,0}=u_c^+\in W_{c,K}$ and $E(z_{c,2})<m^+(c)$, {Lemma~\ref{lem 4.7} implies $z_{c,2} \not \in  W_{c,K}$ and so $\|\nabla z_{c,2}\|_2>t_0$. If $P(z_{c,2})\ge0$, Lemma~\ref{lem 4.5} gives $s^+(z_{c,2})<0\le s^-(z_{c,2})$ and
\[
 E(z_{c,2})\ge E(s^+(z_{c,2})\star z_{c,2})\ge m^+(c),
\]
a contradiction. Thus $P(z_{c,2})<0$, and so Lemma~\ref{lem 4.1} gives $\|\nabla z_{c,2}\|_2 \not =t_0$. At the last $t_*\in(0,2)$ with $\|\nabla z_{c,t_*}\|_2=t_0$, Lemma~\ref{lem 4.1} ensures again that  $P(z_{c,t_*})>0$.} By continuity, some $t_-\in(t_*,2)$ satisfies $P(z_{c,t_-})=0$ and $\|\nabla z_{c,t_-}\|_2>t_0$; hence $z_{c,t_-}\in\mathcal P_c^-$. Therefore
\[
\begin{aligned}
m^-(c)
&\le\max_{t\in[0,2]}E(z_{c,t})\\
&<m^+(c)+\frac{q-1}{2q}S_\beta^{\frac q{q-1}}\\
&=m^+(c)+\frac{\beta+2}{2(N+\beta)}
S_\beta^{\frac{N+\beta}{\beta+2}}.
\end{aligned}
\]
This concludes the proof of the case $N\ge5$.
\end{proof}

\begin{theorem}
\label{thm:second-normalized-solution}
Assume that $\mu>0$, $K>\frac{q}{q-1}$, and
$0<c<\min\{c_F,c_G,c_K\}$. Then $m^-(c)$ is attained
by a positive radial normalized solution $u_c^-\in\mathcal P_c^-$
of problem~\eqref{eq main}, associated with some
$\lambda_c^-\in\mathbb R$. In addition,
\[
 m^+(c)<m^-(c)
 <m^+(c)+\frac{\beta+2}{2(N+\beta)}
 S_\beta^{\frac{N+\beta}{\beta+2}}.
\]
\end{theorem}

\begin{proof}
We first obtain the strict energy gap for the negative branch.
Let $z_{c,t}$, $0\le t\le2$, be the normalized test family
used in the proof of Lemma~\ref{lem 4.10}. The estimates
established there yield, for all sufficiently small
$\varepsilon>0$,
\[
 z_{c,0}=u_c^+,\qquad E(z_{c,2})<m^+(c),\qquad
 \max_{t\in[0,2]}E(z_{c,t})
 <m^+(c)+\frac{q-1}{2q}S_\beta^{\frac q{q-1}}.
\]
By Lemma~\ref{lem 4.10} and for $N \ge 3$, we have
\begin{equation}\label{eq:second-solution-gap}
 m^-(c)\le\max_{t\in[0,2]}E(z_{c,t})
 <m^+(c)+\frac{q-1}{2q}S_\beta^{\frac q{q-1}}.
\end{equation}
Let $\{u_n\}\subset\mathcal P_c^-$ be the sequence from
Lemma~\ref{lem 4.8}. In particular,
\[
\begin{cases}
 E(u_n)\to m^-(c),\\
 P(u_n)=0,\\
 (E|_{S_c})'(u_n)\to0,\\
 Q(u_n)\le-\delta<0.
\end{cases}
\]
The boundedness of $\{u_n\}$ in $H_r$ is part of
Lemma~\ref{lem 4.8}. The Lagrange multiplier argument
used in the proof of Lemma~\ref{lem 4.6} also gives a
bounded sequence $\{\lambda_n\}\subset\mathbb R$ such that
\[
 E'(u_n)+\lambda_n u_n\to0\quad\text{in }H_r^*.
\]
Passing to a subsequence, we have
\[
 \begin{cases}
 u_n\rightharpoonup u\quad\text{in }H_r,\\
 u_n\to u\quad\text{in }L^2(\mathbb R^N),\\
 u_n\to u\quad\text{a.e. in }\mathbb R^N,\\
 \lambda_n\to\lambda\in\mathbb R.
 \end{cases}
\]
The strong $L^2$ convergence follows from
Lemma~\ref{lem 2.1}, and it implies $u\in S_c$.
As in the proof of Lemma~\ref{lem 4.6}, we may pass to
the limit in the Euler--Lagrange equation to conclude that
$u$ is a normalized critical point of $E$ on $S_c$.
In particular, Lemma~\ref{lem 2.4} yields $P(u)=0$.\\

\noindent
It remains to exclude a critical defect. The Hilbert-space
splitting, Lemma~\ref{lem 2.5}, and the nonlocal
Brezis--Lieb lemma (Lemma~\ref{lem 2.2}) imply
\begin{align*}
 \|\nabla u_n\|_2^2
 &=\|\nabla u\|_2^2+\|\nabla(w_n)\|_2^2+o(1),\\
 B_{\alpha,p}[u_n]&=B_{\alpha,p}[u]+o(1),\\
 B_{\beta,q}[u_n]
 &=B_{\beta,q}[u]+B_{\beta,q}[w_n]+o(1).
\end{align*}
Since $P(u_n)=P(u)=0$ and
$\|u_n\|_2^2=\|u\|_2^2=c^2$, it follows that
\begin{equation}\label{eq:second-solution-defect}
 \|\nabla(w_n)\|_2^2
 =B_{\beta,q}[w_n]+o(1),
\end{equation}
where $w_n = u_n - u$. After extraction, we may assume that $\|\nabla(w_n)\|_2^2\to\ell\ge0$ and $B_{\beta,q}[w_n]\to\ell$. By the sharp HLS inequality \eqref{S}, we get
\begin{equation}\label{eq:second-solution-bubble}
 \ell=0\qquad\text{or}\qquad
 \ell\ge S_\beta^{\frac q{q-1}}.
\end{equation}
Testing the approximate Euler--Lagrange equation with $u_n$ and the limit equation with $u$, using the convergence of $\lambda_n$, the lower-order Choquard term and the preceding splittings, we obtain
\[
 \int_{\mathbb R^N}|u_n|^2\log|u_n|^2\,dx
 \to
 \int_{\mathbb R^N}|u|^2\log|u|^2\,dx.
\]
Consequently, the definition of $E$ gives
\begin{equation}\label{eq:second-solution-energy-splitting}
 m^-(c)=E(u)+\frac{q-1}{2q}\ell.
\end{equation}
As $u\in\mathcal P_c$ and $\mathcal P_c^0=\varnothing$,
we have either $u\in\mathcal P_c^+$, in which case
$E(u)\ge m^+(c)$, or $u\in\mathcal P_c^-$, in which
case $E(u)\ge m^-(c)>m^+(c)$ by
Lemma~\ref{lem 4.7}. Thus $E(u)\ge m^+(c)$.
If $\ell>0$, equations
\eqref{eq:second-solution-bubble} and
\eqref{eq:second-solution-energy-splitting} imply
\[
 m^-(c)\ge m^+(c)
 +\frac{q-1}{2q}S_\beta^{\frac q{q-1}},
\]
which contradicts \eqref{eq:second-solution-gap}.
Therefore $\ell=0$ and $u_n\to u$ strongly in
$H^1(\mathbb R^N)$. The convergence of the logarithmic
integrals and of the other energy terms, together with
the Orlicz recovery argument used in the proof of
Lemma~\ref{lem 4.6}, yield $u_n\to u$ strongly in $H_r$.
In particular,
\[
 \begin{cases}
 E(u)=m^-(c),\\
 P(u)=0,\\
 Q(u)=\lim_{n\to\infty}Q(u_n)\le-\delta<0.
 \end{cases}
\]
Hence $u\in\mathcal P_c^-$.\\

\noindent
Finally, since $E(|u|)=E(u)$ and $P(|u|)=P(u)$,
$Q(|u|)=Q(u)$, the function $|u|$ is also a minimizer
of $E$ on $\mathcal P_c^-$. {Since $Q(|u|)<0$, Lemma~\ref{lem 2.4} implies that $|u|$ is a normalized critical point of $E$ on $S_c$.
Therefore, $|u|$ is a nonnegative normalized solution of
\eqref{eq main}.} Applying
Lemma~\ref{lem 2.6} and the same strong maximum
principle argument as in Lemma~\ref{lem 4.6}, we obtain
$|u(x)|>0$ for all $x\in\mathbb R^N$.\\

\noindent
We denote this solution by $u_c^-$ and its Lagrange
multiplier by $\lambda_c^-$. The strict inequality
$m^-(c)>m^+(c)$ follows from Lemma~\ref{lem 4.7}.
\end{proof}

\begin{proof}[Proof of Theorem~\ref{thm 1}:] The existence of $u_c^+$ is established in
Lemma~\ref{lem 4.6}; the existence of $u_c^-$ and the
remaining assertions follow from
Theorem~\ref{thm:second-normalized-solution}.
\end{proof}

\FloatBarrier

\section*{Declarations}

\begin{itemize}
\item{Conflicts/Competing Interests}: The authors declare no conflicts of interest and no competing interests.
\item{Funding}: No external funding was received.
\item{Data Availability}: No data were used in this study.
\item{Author Contributions}: (M.S.A) Methodology and writing--original draft; (R.S. and C.V.) supervision and project administration; (M.B.G. and D.O.) writing--review and editing. All authors have read and approved the final manuscript.
\end{itemize}

\footnotesize{}

\begin{thebibliography}{99}

\bibitem{MSA2026} Abolhassanifar, Mohammad Saeid; Saadati, Reza; Ghaemi, Mohammad Bagher; O'Regan, Donal. Normalized solutions for a critical Schr\"odinger--Choquard equation with double nonlocal terms and a logarithmic nonlinearity. J. Math. Anal. Appl. 562 (2026), no. 1, Paper No. 130707, 43 pp.


\bibitem{Alves2020} Alves, Claudianor O.; Ji, Chao. Existence and concentration of positive solutions for a logarithmic Schrödinger equation via the penalization method. Calc. Var. Partial Differential Equations 59 (2020), no. 1, Paper No. 21, 27 pp.

\bibitem{BartschSoave2017} Bartsch, Thomas; Soave, Nicola. A natural constraint approach to normalized solutions of nonlinear Schr\"odinger equations and systems. J. Funct. Anal. 272 (2017), no. 12, 4998--5037.

\bibitem{BartschLiuLiu2020} Bartsch, Thomas; Liu, Yanyan; Liu, Zeng. Normalized solutions for a class of nonlinear Choquard equations. SN Partial Differential Equations and Applications 1 (2020), Paper No. 34.


\bibitem{Bialynicki1979} Bialynicki-Birula, Iwo; Mycielski, Jerzy. Nonlinear wave mechanics. Phys. Scr. 20 (1979), 539--544.


\bibitem{Cao2025} Cao, Qing-Hai; Yuan, Wen-Shuo; Ge, Bin; Wang, Mei-Yan. On the existence of normalized solutions to the $p$-Laplacian Choquard equation with logarithmic nonlinearity. Anal. Math. Phys. 15 (2025), no. 6, Paper No. 133.

\bibitem{Cazenave1983} Cazenave, Thierry. Stable solutions of the logarithmic Schrödinger equation. Nonlinear Anal. 7 (1983), no. 10, 1127--1140.

\bibitem{ChenJinRadulescuZhou2026} Chen, Sitong; Jin, Peng; R\u{a}dulescu, Vicen\c{t}iu D.; Zhou, Xin'ao. Normalized solutions for the upper critical Choquard equation with nonlocal perturbation. Math. Ann. 395 (2026), Paper No. 20, 52 pp.

\bibitem{ChenZhang2026} Chen, Xingwen; Zhang, Qiongfen. Normalized solutions for a kind of mass--supercritical Schr\"odinger--Choquard equations. AIMS Math. 11 (2026), no. 2, 4656--4680.

\bibitem{Deng2025} Deng, Yinbin; Shi, Yulin; Yang, Xiaolong. Normalized solutions to the Choquard equation, including the critical exponents and a logarithmic perturbation. J. Differential Equations 429 (2025), 204--246.

\bibitem{Deng2021} Deng, Yinbin; Pi, Huirong; Shuai, Wei. Multiple solutions for logarithmic Schrödinger equations with critical growth. Methods Appl. Anal. 28 (2021), no. 2, 221--248.

\bibitem{DingWang2024} Ding, Yanheng; Wang, Hua-Yang. Normalized solutions to Schr\"odinger equations with critical exponent and mixed nonlocal nonlinearities. J. Geom. Anal. 34 (2024), no. 7, Paper No. 215.

\bibitem{GaoHe2024} Gao, Qian; He, Xiaoming. Normalized solutions for the Choquard equations with critical nonlinearities. Adv. Nonlinear Anal. 13 (2024), no. 1, Paper No. 20240030.

\bibitem{Hongyu2016} Ye, Hongyu. Mass minimizers and concentration for nonlinear Choquard equations in $\mathbb{R}^N$. Topol. Methods Nonlinear Anal. 48 (2016), no. 2, 393--417.

\bibitem{Lieb2001} Lieb, Elliott H.; Loss, Michael. Analysis. Second edition. Graduate Studies in Mathematics, 14. American Mathematical Society, Providence, RI, 2001. xxii+346 pp.

\bibitem{Lieb1983} Lieb, Elliott H. Sharp constants in the Hardy-Littlewood-Sobolev and related inequalities. Ann. of Math. (2) 118 (1983), no. 2, 349--374.

\bibitem{LiuZhangSun2025} Liu, Jianlun; Zhang, Ziheng; Sun, Hong-Rui. Normalized ground state solutions to HLS upper critical Choquard equation with mass critical or supercritical and Sobolev critical growth. Math. Methods Appl. Sci. 48 (2025), no. 15, 14276--14289.

\bibitem{MaMoroz2024} Ma, Shiwang; Moroz, Vitaly. Asymptotic profiles for Choquard equations with combined attractive nonlinearities. J. Differential Equations 412 (2024), 613--689.

\bibitem{MaMoroz2025} Ma, Shiwang; Moroz, Vitaly. Asymptotic profiles for Choquard equations with combined attractive nonlinearities: the locally critical case. Calc. Var. Partial Differential Equations 64 (2025), Paper No. 127.

\bibitem{MorozVanSchaftingen2015} Moroz, Vitaly; Van Schaftingen, Jean. Existence of ground states for a class of nonlinear Choquard equations. Trans. Amer. Math. Soc. 367 (2015), no. 9, 6557--6579.

\bibitem{Moroz2013} Moroz, Vitaly; Van Schaftingen, Jean. Ground states of nonlinear Choquard equations: existence, qualitative properties and decay asymptotics. J. Funct. Anal. 265 (2013), no. 2, 153--184. 

\bibitem{Moroz2017} Moroz, Vitaly; Van Schaftingen, Jean. A guide to the Choquard equation. J. Fixed Point Theory Appl. 19 (2017), no. 1, 773--813.

\bibitem{Pekar1946} Pekar, S. I. Untersuchung über die Elektronentheorie der Kristalle. J. Phys. USSR 10 (1946), 341.

\bibitem{ShangMa2023} Shang, Xudong; Ma, Pei. Normalized solutions to the nonlinear Choquard equations with Hardy--Littlewood--Sobolev upper critical exponent. J. Math. Anal. Appl. 521 (2023), no. 2, Paper No. 126916.

\bibitem{Shen2025} Shen, Liejun; Squassina, Marco. Existence and concentration of normalized solutions for $p$-Laplacian equations with logarithmic nonlinearity. J. Differential Equations 421 (2025), 1--49.

\bibitem{Shuai2025} Shuai, Wei; Yang, Xiaolong. Normalized solutions for logarithmic Schrödinger equation with a perturbation of power law nonlinearity. Ann. Fenn. Math. 50 (2025), no. 1, 301--330.

\bibitem{Vazquez1984} V\'azquez, Juan Luis. A strong maximum principle for some quasilinear elliptic equations. Appl. Math. Optim. 12 (1984), 191--202.

\bibitem{WangMaLiu2023} Wang, Yachen; Ma, Shiwang; Liu, Xiaonan. Asymptotic behaviors of normalized solutions for a class of Choquard equations. Appl. Math. Lett. 142 (2023), Paper No. 108638.

\bibitem{Weinstein1982} Weinstein, Michael I. Nonlinear Schrödinger equations and sharp interpolation estimates. Comm. Math. Phys. 87 (1982/83), no. 4, 567--576.

\bibitem{XiaZhang2023} Xia, Jiankang; Zhang, Xu. Normalized saddle solutions for a mass supercritical Choquard equation. J. Differential Equations 364 (2023), 471--497.

\bibitem{YangYang2025} Yang, Kaizhen; Yang, Xianyong. Normalized ground state solutions for a mass supercritical Choquard equation. Bull. Math. Sci. 15 (2025), Paper No. 2550022.

\bibitem{YaoChenRadulescuSun2022} Yao, Shuai; Chen, Haibo; R\u{a}dulescu, Vicen\c{t}iu D.; Sun, Juntao. Normalized solutions for lower critical Choquard equations with critical Sobolev perturbation. SIAM J. Math. Anal. 54 (2022), no. 3, 3696--3723.

\bibitem{Ye2022} Ye, Weiwei; Shen, Zifei; Yang, Minbo. Normalized solutions for a critical Hartree equation with perturbation. J. Geom. Anal. 32 (2022), no. 9, Paper No. 242, 44 pp.

\end{thebibliography}
\end{document}